\documentclass[12pt,pdftex]{articleplusnothms}
\usepackage[utf8]{inputenc}
\author{Ileana Streinu and Louis Theran}
\title{
Sparse Hypergraphs and Pebble Game Algorithms}

\usepackage{subfigure}
\withcomplaints
\newcommand{\peb}{\ensuremath{\operatorname{peb}}}
\newcommand{\grsp}{\ensuremath{\operatorname{span}}}
\newcommand{\out}{\ensuremath{\operatorname{out}}}

\newtheorem{theorem}{Theorem}[section]
\newtheorem{corollary}[theorem]{Corollary}
\newtheorem{lemma}[theorem]{Lemma}
\newtheorem{proposition}[theorem]{Proposition}
\newtheorem{algorithm}[theorem]{Algorithm}
\newcommand{\labelprop}[1]{\label{prop.#1}}
\newcommand{\refprop}[1]{Proposition \ref{prop.#1}}

\newcommand{\restateenv}{ZZZ}
\newenvironment{restate}[2]{
  \renewcommand{\restateenv}{restate.#1}
  \newtheorem*{\restateenv}{\refthm{#1}}
  \begin{\restateenv}[#2]
}{\end{\restateenv}}

\date{}

\begin{document}
\maketitle

\begin{abstract}
A hypergraph $G=(V,E)$ is $(k,\ell)$-sparse if no subset $V'\subset
V$ spans more than $k|V'|-\ell$ hyperedges. We characterize
$(k,\ell)$-sparse hypergraphs in terms of graph theoretic, matroidal
and algorithmic properties. We extend several well-known theorems of
Haas, Lov{\'{a}}sz, Nash-Williams, Tutte, and White and Whiteley,
linking arboricity of graphs to certain counts on the number of
edges. We also address the problem of finding lower-dimensional
representations of sparse hypergraphs, and identify a critical
behaviour in terms of the sparsity parameters $k$ and $\ell$. Our
constructions extend the pebble games of Lee and Streinu
\cite{LeSt05} from graphs to hypergraphs.
\end{abstract}

\section{Introduction \labelsec{introduction}}

The focus of this paper is on $(k,\ell)$-sparse hypergraphs.  A
hypergraph (or set system) is a pair $G=(V,E)$ with {\bf vertices}
$V$, $n=|V|$ and {\bf edges} $E$ which are subsets of $V$ (multiple
edges are allowed). If all the edges have exactly two vertices, $G$
is a (multi){\bf graph}. We say that a hypergraph is
$(k,\ell)$-sparse if no subset $V'\subset V$ of $n'=|V'|$ vertices
spans more than $kn'-\ell$ edges in the hypergraph. If, in addition,
$G$ has exactly $kn-\ell$ edges, we say it is $(k,\ell)$-tight.

The $(k,\ell)$-sparse graphs and hypergraphs have applications in
determining connectivity and arboricity (defined later). For some
special values of $k$ and $\ell$, the $(k,\ell)$-sparse graphs have
important applications to rigidity theory: bar-and-joint minimally
rigid frameworks in dimension 2, and body-and-bar structures in
arbitrary dimension are both characterized generically by sparse
graphs.

In this paper, we prove several equivalent characterizations of the
$(k,\ell)$-sparse hypergraphs, and give efficient algorithms for
three specific problems.  The {\bf decision} problem asks if a
hypergraph $G$ is $(k,\ell)$-tight.  The {\bf extraction} problem
takes an arbitrary hypergraph $G$ as input and returns as output a
maximum size (in terms of edges) $(k,\ell)$-sparse sub-hypergraph of
$G$. The {\bf components} problem takes a {\em sparse} $G$ as input
and returns as output the {\em maximal} $(k,\ell)$-tight induced
sub-hypergraphs of $G$.

The {\em dimension} of a hypergraph is its minimum edge size. A
large dimension makes them difficult to visualize. We also address
the {\bf representation} problem, which asks for finding a suitably
defined lower-dimensional hypergraph in the same sparsity class, and
we identify a critical behaviour in terms of the sparsity parameters
$k$ and $\ell$.

There is a vast literature on sparse $2$-graphs (see Section
\ref{sec.related}), but not so much on hypergraphs. In this paper,
we carry over to the most general setting the characterization of
sparsity via pebble games from Lee and Streinu \cite{LeSt05}. Along
the way, we develop structural properties for sparse hypergraph
decompositions, identify the problem of lower dimensional
representations, give the proper hypergraph version of depth-first
search in a directed sense and  apply the pebble game to efficiently
find lower-dimensional representations within the same sparsity
class.

Complete historical background is given in Section \ref{sec.related}.
In Section 2, we describe our pebble game for hypergraphs
in detail.  The rest of the paper provides the proofs: Sections 3
and 4 address structural properties of sparse hypergraphs; Sections 5 and
6 relate graphs accepted by the pebble game with sparse hypergraphs; Section 7
addresses the questions of representing sparse hypergraphs by lower dimensional
ones.

\subsection{Preliminaries and related work\labelsec{preliminaries}}
In this section we give the definitions and describe the notation
used in the paper.

{\bf Note:} for simplification, we will often use {\em graph}
instead of {\em hypergraph} and {\em edge} instead of {\em
hyperedge}, when the context is clear.

\paragraph{Hypergraphs. \labelsec{hypergraphs}} Let $G=(V,E)$ be a
hypergraph, i.e. the edges of $G$ are subsets of $V$. A vertex $v\in
e$ is called an {\em endpoint} (or simply {\em end}) of the edge. We
allow parallel edges, i.e. multiple copies of the same edge.

\begin{figure}%[htbp]
\centering
\subfigure[]{\includegraphics[height=1.0 in]{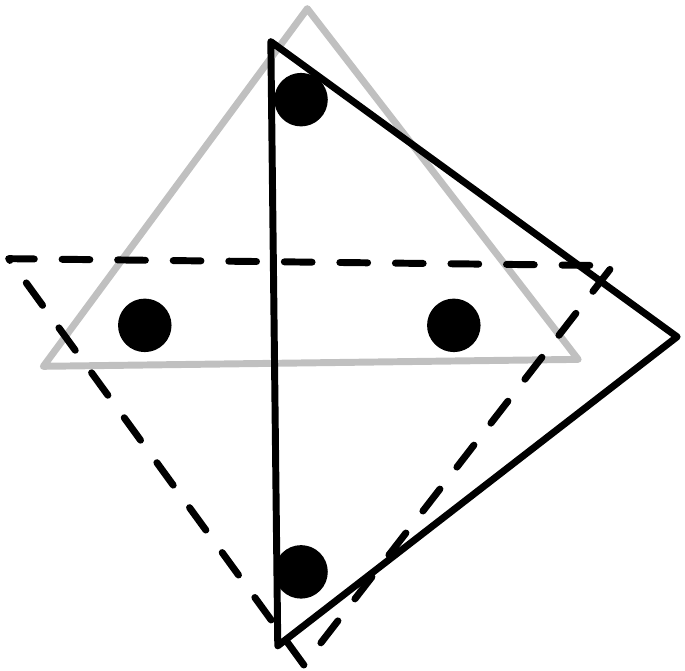}}
\hspace{.3in}
\subfigure[]{\includegraphics[height=1.0 in]{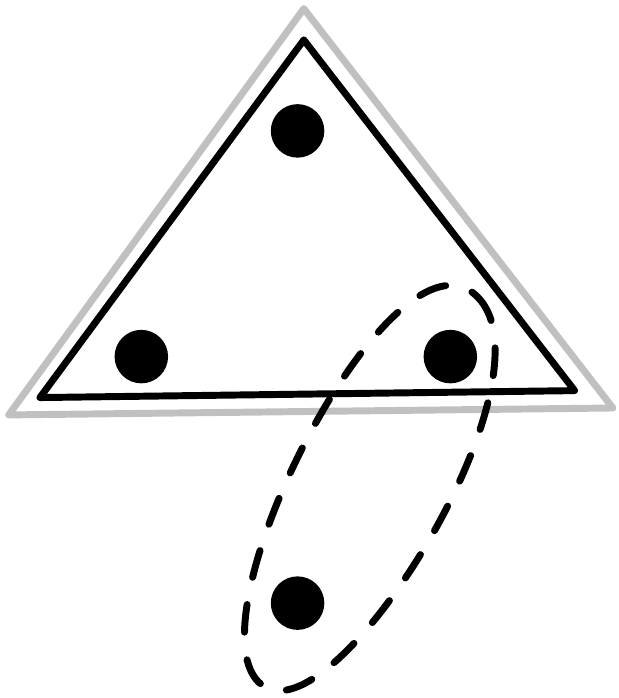}}
\caption{Two hypergraphs.  The hypergraph in (a) is 3-uniform; (b) is 2-dimensional but not
a 2-graph.}
\label{fig.hypergraph-examples}
\end{figure}

For a subset $V'$ of the vertex set $V$, we define span($V'$), the
{\bf span} of $V'$, as the set of edges with endpoints in $V'$:
$E(V')=\{e\in E : e\subset V'\}$. Similarly, for a subset $E'$ of
$E$, we define the span of $E'$ as the set of vertices in the union
of the edges: $V(E')=\bigcup_{e\in E'} e$.  The {\bf hypergraph
dimension} (or dimension) of an edge is its number of elements. The
hypergraph dimension of a graph $G$ is its {\em minimum} edge
dimension. A graph in which each edge has dimension $s$ is called
{\bf $s$-uniform} or, more succinctly, a {\bf $s$-graph}. So what is
typically called a graph in the literature is a $2$-graph, in our
terminology. \reffig{hypergraph-examples} shows two examples of
hypergraphs.

\begin{figure}%[htbp]
\centering %%
\subfigure[]{\label{fig.represented-example-1}\includegraphics[height=1.0
in]{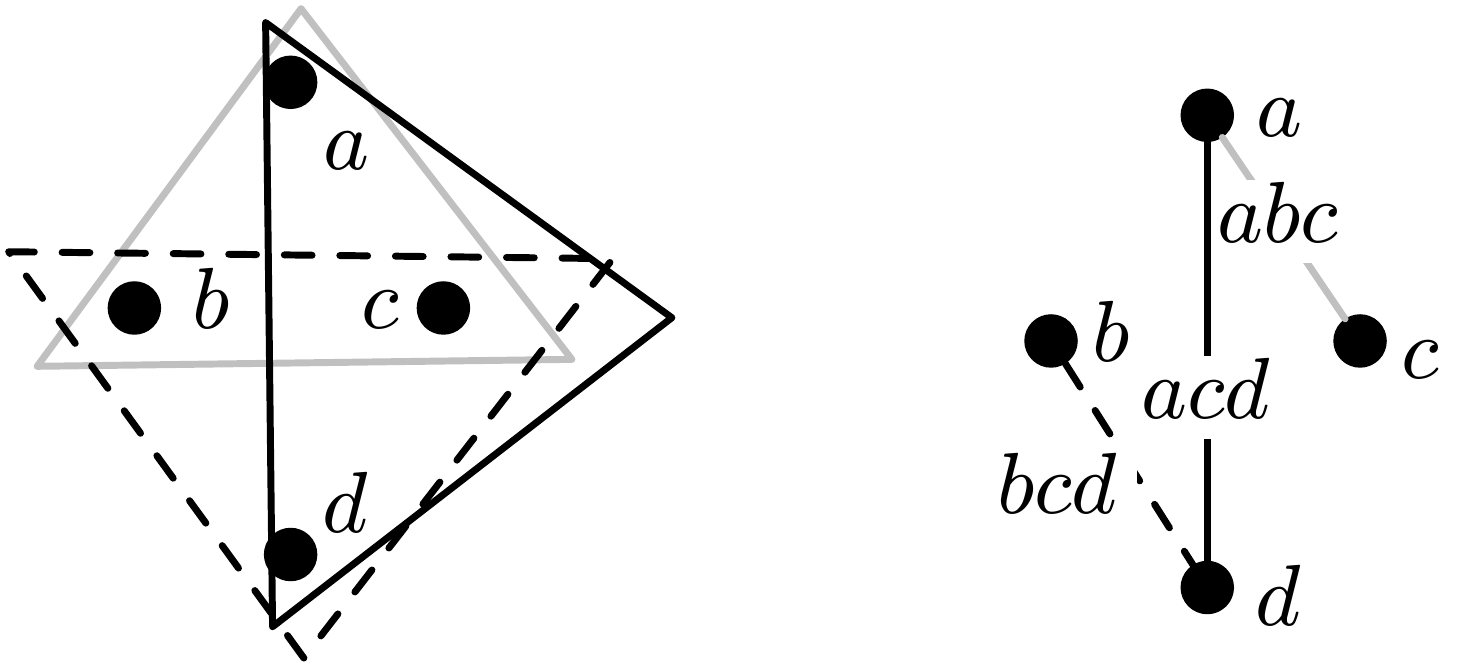}}
\hspace{.3in}
\subfigure[]{\label{fig.represented-example-2}\includegraphics[height=1.0
in]{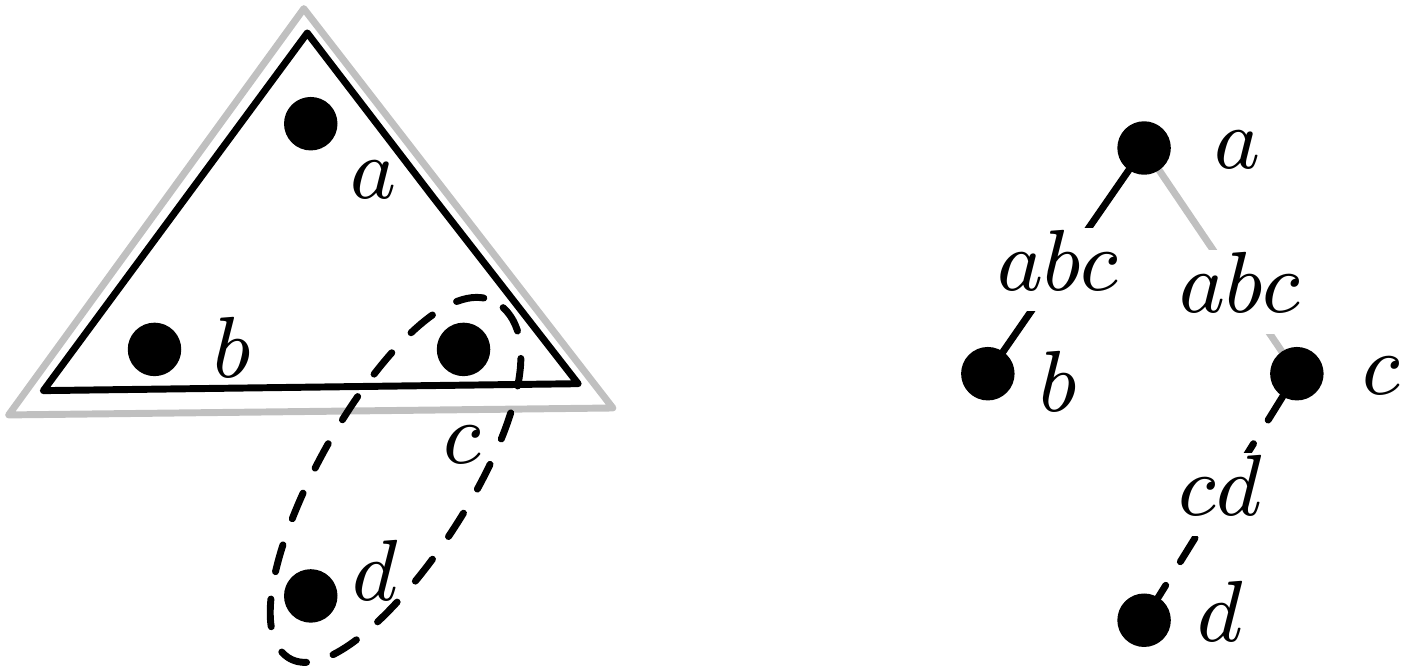}} \caption{Lower dimensional
representations.  In both cases, the $2$-uniform graph  on the right
(a tree) represents the hypergraph on the left (a hypergraph tree)
with respect to $(1,1)$-sparsity. The $2$-dimensional
representations of edges have similar styles to the edges they represent
 and are labeled with the vertices of the hyperedge.}
\label{fig.representations}
\end{figure}

\begin{figure}%[htbp]
\centering
\includegraphics[height=1 in]{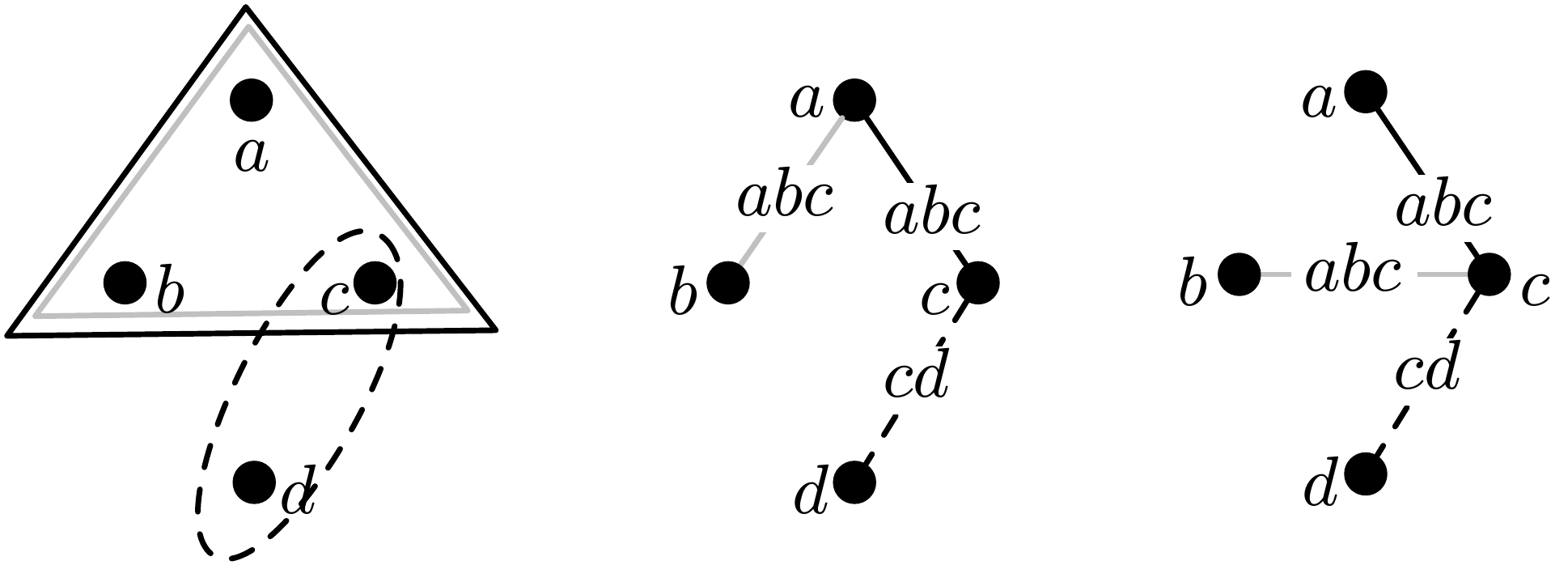}
\caption{Lower dimensional representations are not unique.  Here we
show two 2-uniform representations of the same hypergraph with
respect to $(1,1$)-sparsity.}\label{fig.notunique}
\end{figure}

We say that a hypergraph $H=(V,F)$ {\bf represents} a hypergraph
$G=(V,E)$ with respect to some property ${\cal P}$, if both $H$ and
$G$ satisfy the property, and there is an isomorphism $f$ from $E$
to $F$ such that $f(e)\subset e$ for all $e\in E$. In this paper, we
are primarily concerned with representations which preserve
sparsity. In our figures, we visually present hypergraphs as their
lower dimensional representations when possible, as in
\reffig{representations}.  We observe that representations with
respect to sparsity are not unique, as shown in \reffig{notunique}.

The standard concept of {\bf degree} of a vertex $v$ extends
naturally to hypergraphs, and is defined as the number of edges to
which $v$ belongs. The degree of a set of vertices $V'$ is the
number of edges with at least one endpoint in $V'$ and another in
$V-V'$.

\begin{figure}%[htbp]
\centering %%
\includegraphics[height=1.0 in]{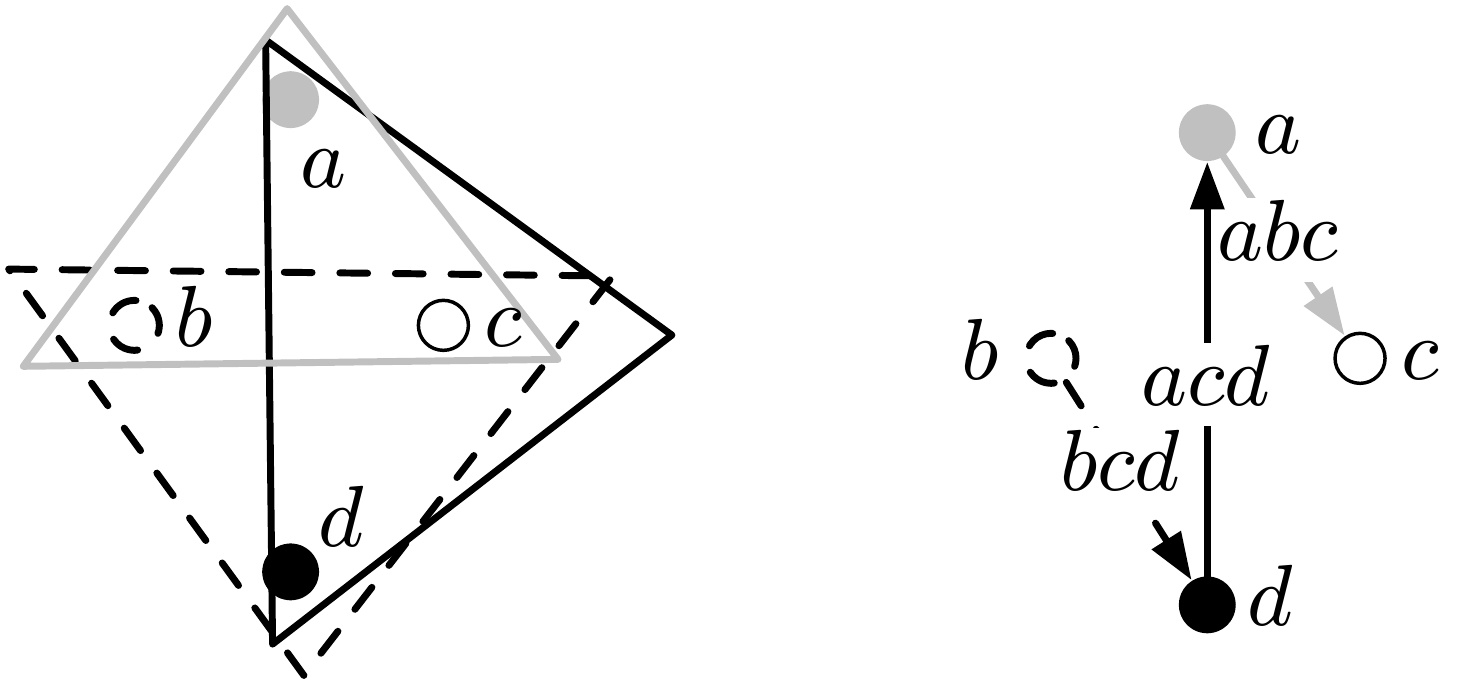}
\caption{An oriented 3-uniform hypergraph. On the left, the tail of each edge is indicated by the style of the vertex.
In the 2-uniform representation on the right, the edges are shown as directed arcs.}
\label{fig.oriented-example}
\end{figure}

An {\bf orientation} of a hypergraph is given by identifying as the
{\bf tail} of each edge one of its endpoints.
\reffig{oriented-example} shows an oriented hypergraph and a lower
dimensional representation of the same graph.

In an oriented hypergraph, a {\bf path} from a vertex $v_1$ to a
vertex $v_t$ is given by a sequence
\begin{eqnarray}
v_1,e_1,v_2,e_2,\ldots,v_{t-1},e_{t-1},v_t
\end{eqnarray}
where $v_i$ is an endpoint of $e_{i-1}$ and $v_i$ is the tail of
$e_i$ for $1\le i\le t-1$.

The concepts of in-degree and out-degree extend to oriented
hypergraphs. The out-degree of a vertex is the number of edges which
identify it as the tail and connect $v$ to $V-v$; the in-degree is
the number of edges that do not identify it as the tail. The
out-degree of a subset $V'$ of $V$ is the number of edges with the
tail in $V'$ and at least one endpoint in $V-V'$; the in-degree of
$V'$ is defined symmetrically.  It is easy to check that the
out-degree and in-degree of $V'$ sum to the undirected degree of
$V'$. Notice that loops (one-dimensional edges) contribute nothing
to the out-degree of a vertex-set.

We use the notation $N_{G}(V')$ to denote the set of neighbors in
$G$ of a subset $V'$ of $V$.

\begin{figure}[htbp]
\centering %%
\subfigure[]{\label{fig.dfs-1}\includegraphics[height=1 in]{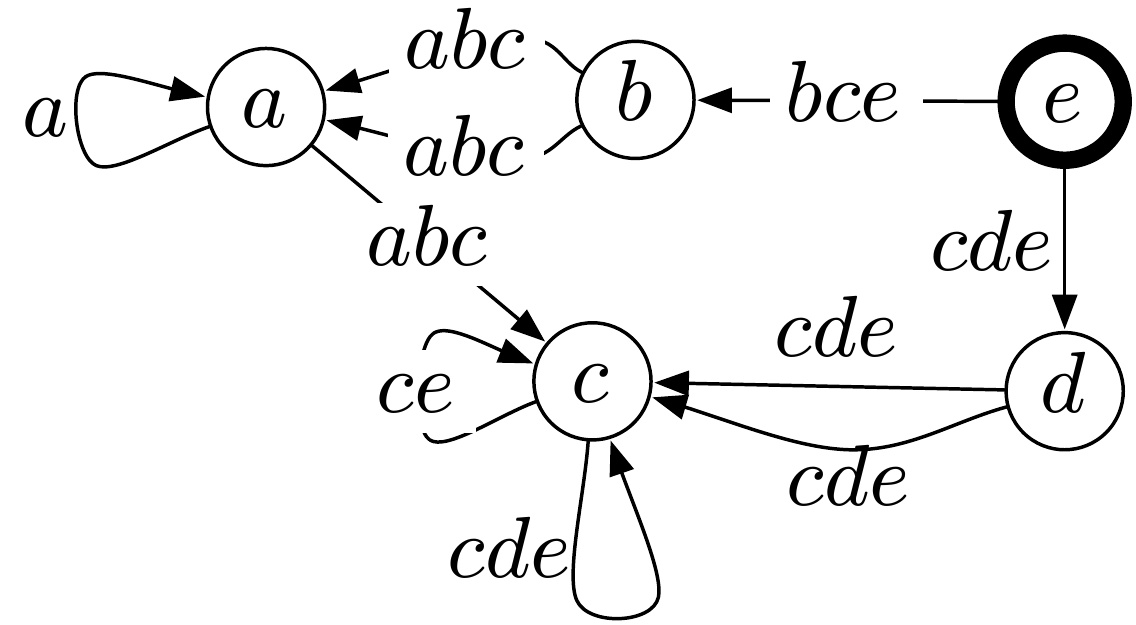}}
%\hspace{.3in}
\subfigure[]{\label{fig.dfs-2}\includegraphics[height=1 in]{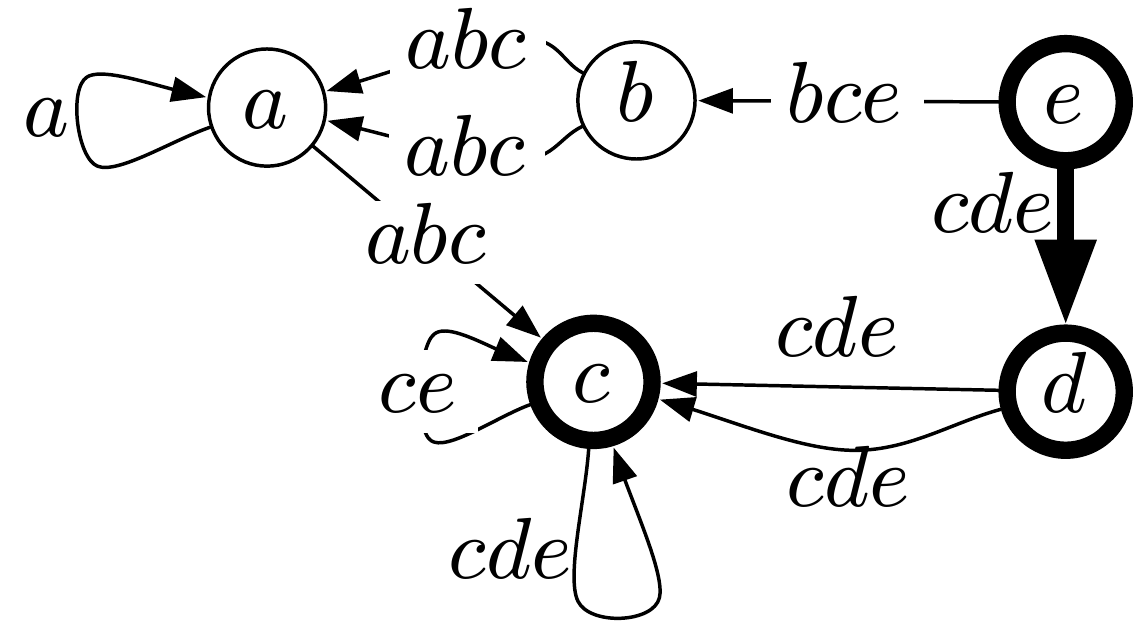}}
%\hspace{.3in}
\subfigure[]{\label{fig.dfs-3}\includegraphics[height=1 in]{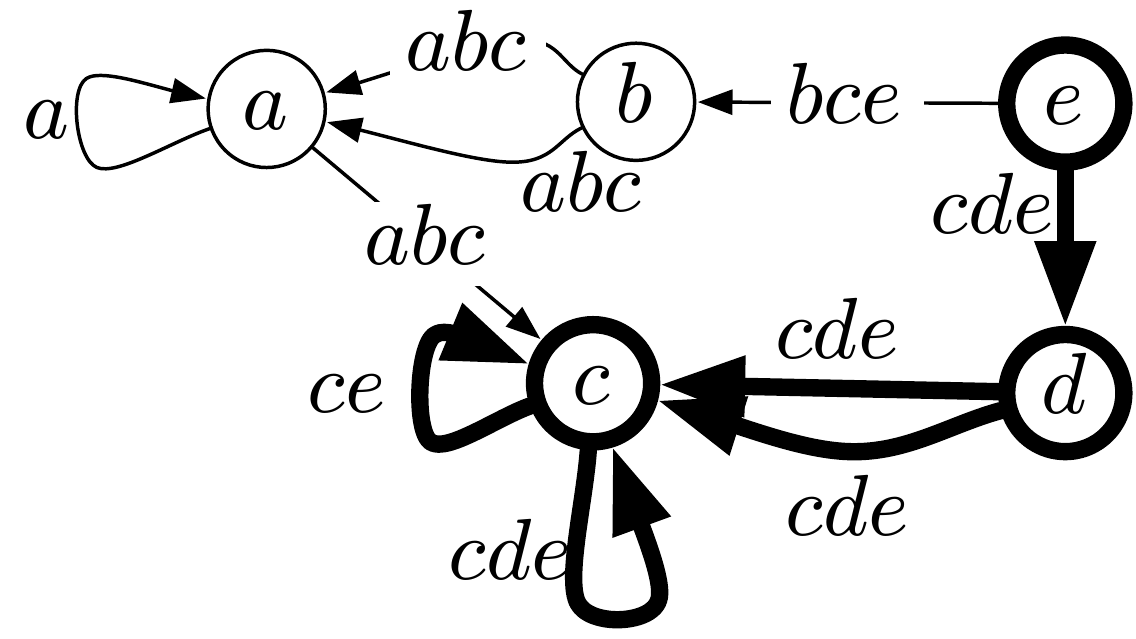}}
%%
%\vspace{.3 in}
\subfigure[]{\label{fig.dfs-4}\includegraphics[height=1 in]{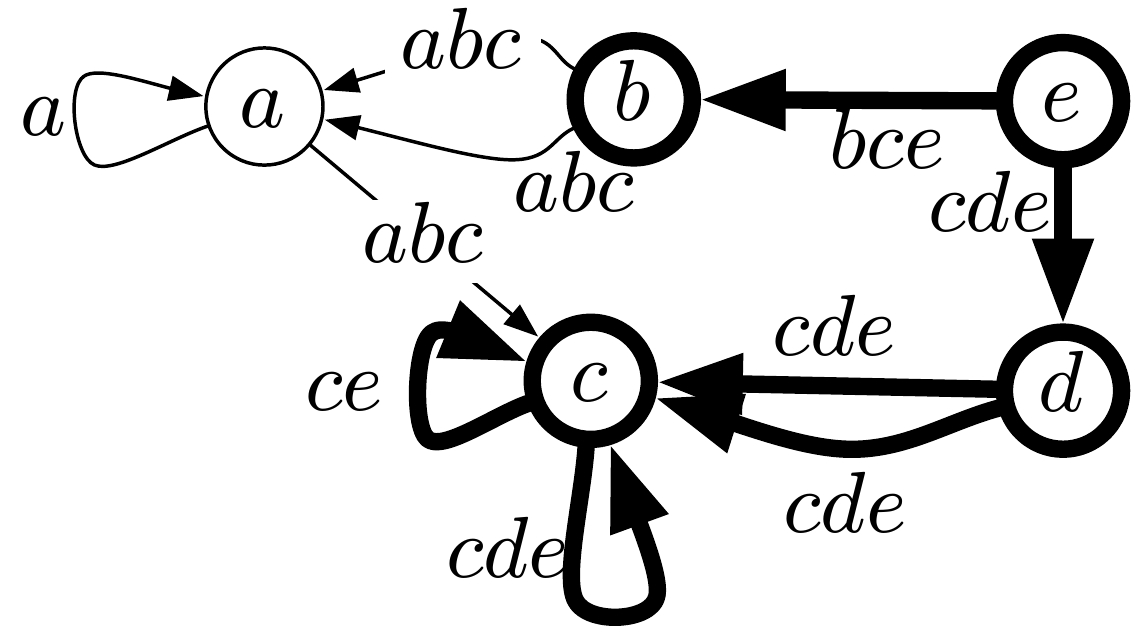}}
%\hspace{.3in}
\subfigure[]{\label{fig.dfs-5}\includegraphics[height=1 in]{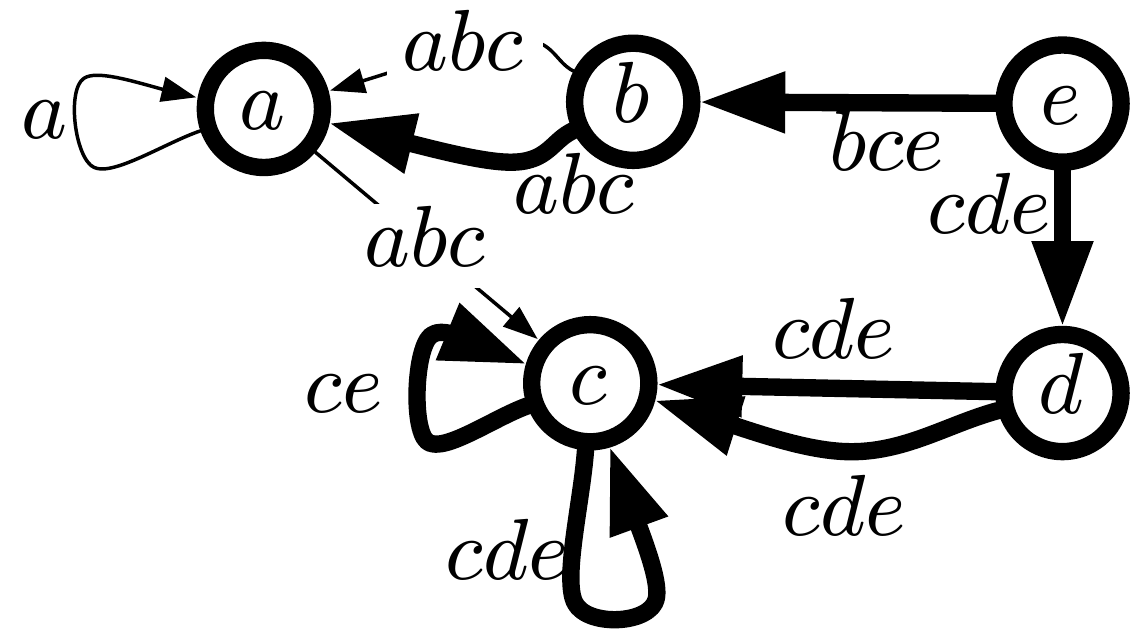}}
%\hspace{.3in}
\subfigure[]{\label{fig.dfs-6}\includegraphics[height=1 in]{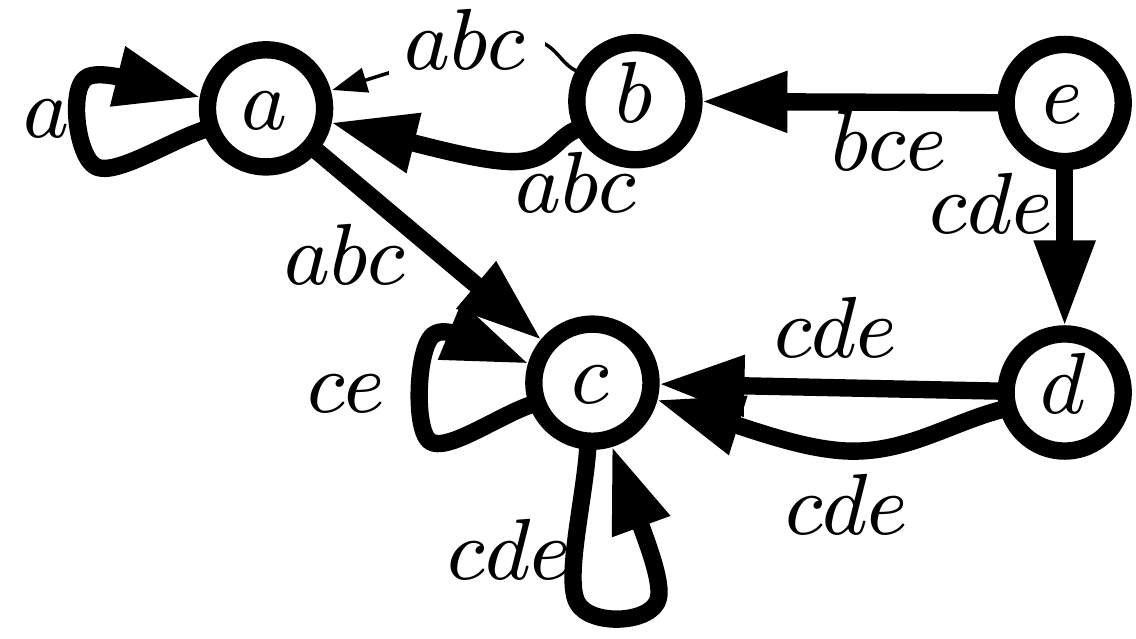}}
\subfigure[]{\label{fig.dfs-7}\includegraphics[height=1 in]{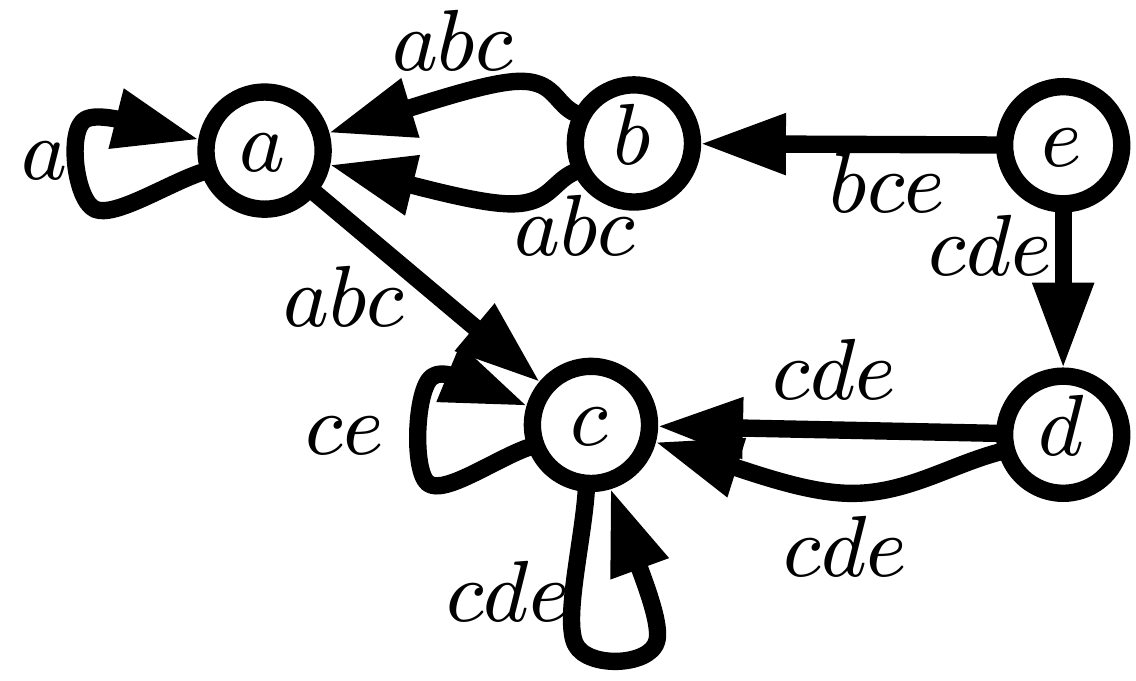}}
\caption{Searching a hypergraph with depth-first search starting at vertex $e$.
Visited edges and vertices are shown with thicker lines.
The search proceeds across an edge from the tail to each of the other endpoints and
backs up at an edge when all its endpoints have been visited
(as in the transition from (b) to (c)). }
\label{fig.dfs}
\end{figure}

The standard depth-first search algorithm in directed graphs,
starting from a source vertex $v$, extends naturally to oriented
hypergraphs: recursively explore the graph from the unexplored
neighbors of $v$, one after another (ending when it has no
unexplored neighbors left). We will use it in the implementation of
the pebble game to explore vertices of hypergraphs.
\reffig{dfs} shows the depth-first exploration of a hypergraph.
Notice that the picture uses a uniform $2$-dimensional
representation for a $3$-hypergraph (the hyperedges should be clear
from the labels on the $2$-edges representing them).

Table \ref{tab.hypergraph-terminology} gives a summary of the
terminology in this section.

\begin{table}
\centering
\begin{tabular}{|l|l|l|}
\hline
{\bf Term} & {\bf Notation} & {\bf Meaning} \\
\hline
\hline
Edge & $e$ & $e\subset V$
\\
\hline
Graph & $G=(V,E)$ &  $V$ is a finite set of vertices ; $E\subset 2^V$ is a set of edges \\
\hline
Subset of vertices & $V'$ & $V'\subset V$ \\
\hline
Size of $V'$ & $n'$ & \card{V'} \\
\hline
Subset of edges & $E'$ & $E'\subset E$ \\
\hline
Size of a subset of edges & $m'$ & $\card{E'}$ \\
\hline
Span of $V'$ & $E(V')$ & Edges in $E$ that are subsets of $V'$ \\
\hline
Span of $E'$ & $V(E')$  & Vertices in the union of $e\in E'$ \\
\hline
Dimension of  $e\in E$ & $|e|$ & Number of elements in $e$ \\
\hline
Dimension of $G$ & $s$  & Minimum dimension of an edge in $E$. \\
\hline
Max size of an edge & $s^*$ & Maximum size of an edge in $E$ \\
\hline
Neighbors of $V'$ in $G$ & $N_G(V')$ & Vertices connected to some $v\in V'$ \\
\hline
\end{tabular}
\caption{Hypergraph terminology used in this paper.}
\label{tab.hypergraph-terminology}
\end{table}

\paragraph{Sparse hypergraphs.\labelsec{sparse}} A graph is {\bf
$(k,\ell)$-sparse} if for any subset $V'$ of $n'$ vertices and its
span $E'$, $m'=|E'|$:
\begin{eqnarray}
m' \le kn'-\ell
\labeleq{subset}
\end{eqnarray}

\begin{figure}%[htbp]
\centering %%
\includegraphics[height=1.5 in]{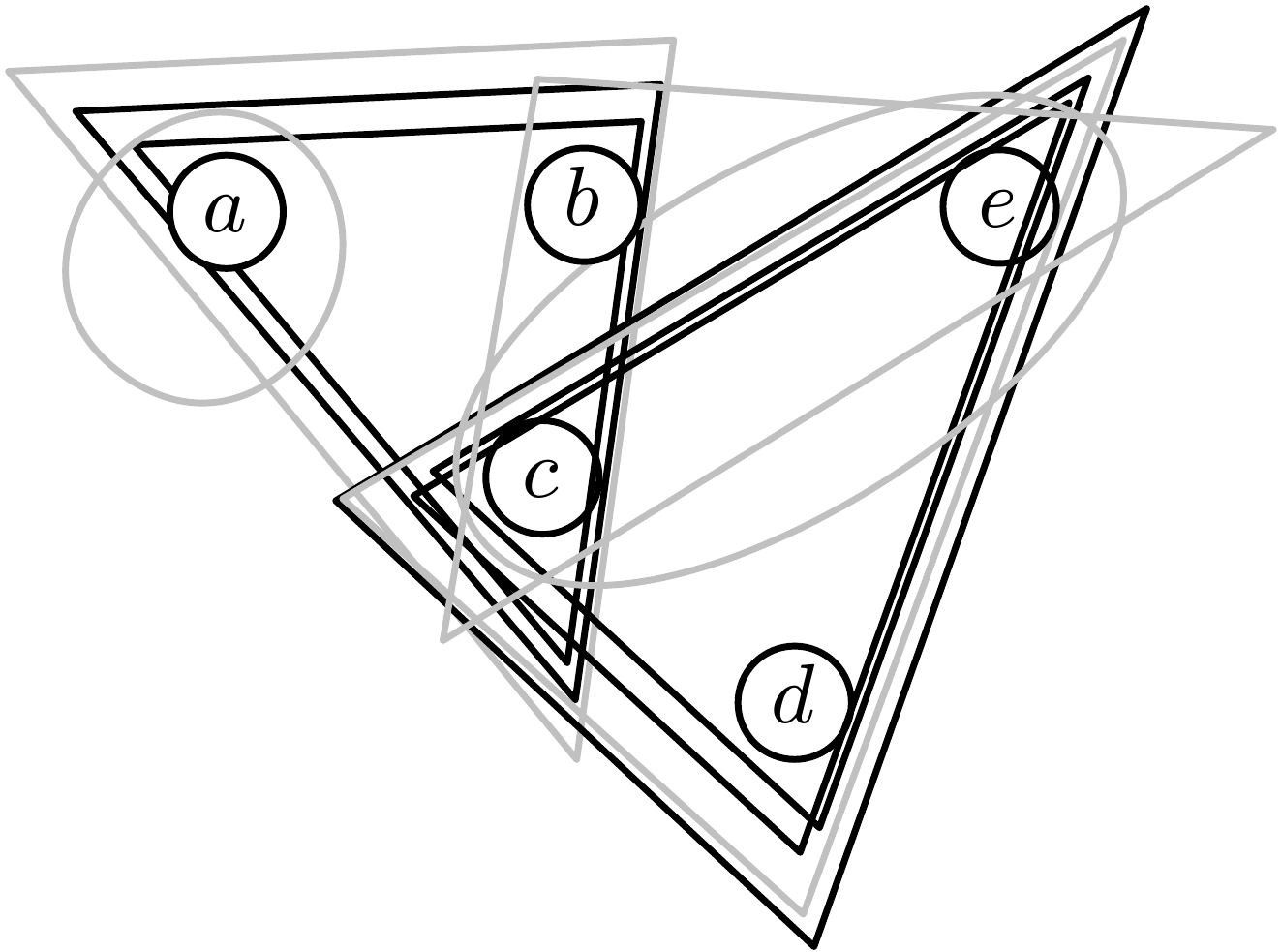}
\caption{A (2,0)-tight hypergraph decomposed into two $(1,0)$-tight ones (gray and black).}
\label{fig.2-map-tight}
\end{figure}

A sparse graph that has exactly $kn-\ell$ edges is called {\bf
tight}; \reffig{2-map-tight} shows a $(2,0)$-tight hypergraph.  A
graph that is not sparse is called {\bf dependent}.

A simple observation, formalized below in
\reflem{sparse-graph-rank}, implies that $0\leq \ell \leq s k -1$,
for sparse hypergraphs of dimension $s$. {\em From now on, we will
work with parameters $k, \ell$ and $s$ satisfying this condition.}

We also define $K_n^{k,\ell}$ as the complete hypergraph with edge
multiplicity $ks-\ell$ for $s$-edges.  For example $K_n^{k,0}$ has:
$k$ loops on every vertex, $2k$ copies of every $2$-edge, $3k$
copies of every $3$-edge, and so on.
\reflem{loops-and-parallel-edges} shows that every sparse graph is a
subgraph of $K_n^{k,0}$.

A sparse graph $G$ is {\bf critical} if the only representation of $G$ that is
sparse is $G$ itself.  In terms of $B_{G}$ this means that no proper
subgraph of $B'$ of $B_{G}$ corresponds to a hypergraph that is
sparse.

There are two important types of subgraphs of sparse graphs.  A {\bf block} is a
tight subgraph of a sparse graph.  A {\bf component}
 is a maximal block.

In this paper, we study five computational problems.  The {\bf decision} problem
asks if a graph $G$ is $(k,\ell)$-tight.  The {\bf extraction} problem takes a graph
$G$ as input and returns as output a maximum $(k,\ell)$-sparse subgraph of
$G$.  The {\bf optimization} problem is a variant of the {\bf extraction} problem;
it takes as its input
a graph $G$ and a weight function on $E$ and returns as
its output a minimum weight maximum $(k,\ell)$-sparse subgraph of $G$.
The {\bf components} problem take a graph $G$ as input and returns
as output the components of $G$.  The {\bf representation} problem takes
as input a sparse graph $G$ and returns as output a sparse graph $H$
that represents $G$ and has lower dimension if this is possible.

\begin{table}
\centering
\begin{tabular}{|l|l|}
\hline
{\bf Term} & {\bf Meaning} \\
\hline
\hline
Sparse graph $G$ & $m'\leq kn'-l$ for all subsets $E'$, $m'=|E'|$. \\
\hline
Tight graph $G$ & $G$ is  sparse with $kn-\ell$ edges. \\
\hline
Dependent graph $G$  & $G$ is not sparse \\
\hline
Block $H$ in $G$ & $G$ is sparse, and $H$ is a tight subgraph \\
\hline
Component $H$ of $G$ & $G$ is sparse and $H$ is a maximal block \\
\hline
Decision problem  & Decide if a graph $G$ is sparse \\
\hline
Extraction problem  & Given $G$, find a maximum sized sparse subgraph $H$ \\
\hline
Optimization problem & Given $G$, find a minimum weight maximum sized sparse subgraph $H$ \\
\hline
Components problem  & Given $G$, find the components of $G$ \\
\hline
Representation problem & Given a sparse $G$, find a sparse representation of lower
dimension\\
\hline
\end{tabular}
\caption{Sparse graph terminology used in this paper.}
\label{tab.sparse-terminology}
\end{table}

Table \ref{tab.sparse-terminology} summarizes the notation and terminology
related to sparseness used in this paper.

While the definitions in this section are made for families of
sparse graphs, they can be interpreted in terms of matroids and
rigidity theory.  Table \ref{tab.sparse-concepts} relates the
concepts in this section to matroids and generic rigidity, and can
be skipped by readers who are not familiar with these fields.

\begin{table}
\centering
\begin{tabular}{|l|l|l|}
\hline
{\bf Sparse graphs} & {\bf Matroids} & {\bf Rigidity}  \\
\hline
\hline
Sparse & Independent & No over-constraints \\
\hline
 Tight & Independent and spanning & Isostatic/minimally rigid \\
\hline
 Block & --- & Isostatic region \\
\hline
Component & --- & Maximal isostatic region \\
\hline
Dependent & Contains a circuit & Has stressed regions \\
\hline
\end{tabular}
\caption{Sparse graph concepts and analogs in matroids and rigidity.}
\label{tab.sparse-concepts}
\end{table}

\begin{figure}[htbp]
\centering %%
\includegraphics[height=1.5 in]{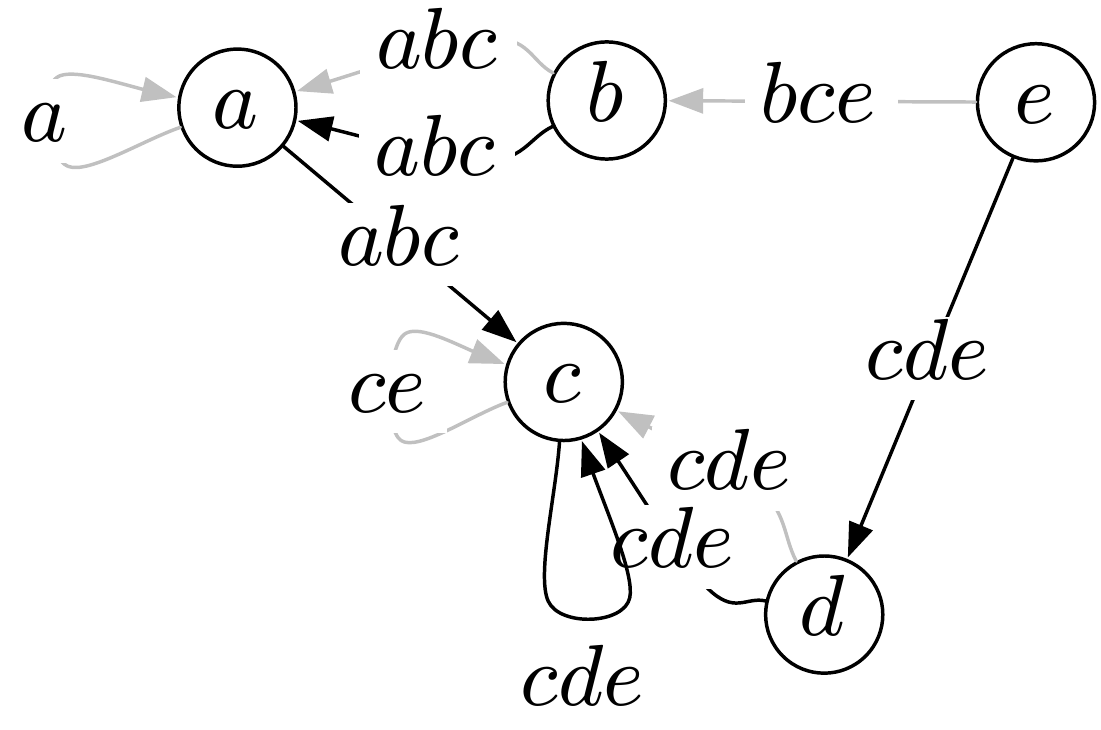}
\caption{The hypergraph from \reffig{2-map-tight}, shown here in a
lower-dimensional representation, is a 2-map.  The maps are black
and gray.
 Observe that each vertex is the tail of one black edge and one gray one.}\label{fig.2-map-oriented}
\end{figure}

\paragraph{Fundamental hypergraphs.} A {\bf map} is a hypergraph that admits
an orientation such that the out degree of every vertex is exactly
one. A $k$-{\bf map} is a graph that admits a decomposition into $k$
disjoint maps.  \reffig{2-map-oriented} shows a $2$-map, with an
 orientation of the edges certifying that the graph is a $2$-map.

An edge $e$ {\bf connects} subsets $X$ and $Y$ of $V$ if $e$ has an
end in both $X$ and $Y$.  A graph is {\bf $k$-edge connected} if
$\card{E(X,V-X)}\ge k,$
for any subset $X$ of $V$, where $E(X,Y)$ is the set of edges connecting
$X$ and $Y$.

A graph is {\bf $k$-partition connected} if
\begin{eqnarray}
\card{\bigcup_{i\neq j} E(P_i,P_j)}\ge k(t-1)
\label{partition}
\end{eqnarray}
for any partition $\mathcal{P}=\{P_1,P_2,\ldots,P_t\}$ of $V$.
This definition appears in  \cite{frank2001}.

A {\bf tree} is a minimally 1-partition connected graph. A reminder
that this is the definition of a {\em tree} in a {\em hypergraph},
but we use the shortened terminology and drop {\em hyper}. A
$k$-arborescence is a graph that admits a decomposition into $k$
disjoint trees. For $2$-graphs, the definitions of partition
connectivity and edge connectivity coincide by the well-known
theorems of Tutte \cite{tutte61} and Nash-Williams \cite{Na61}.  We
also observe that for general hypergraphs, connectivity and
$1$-partition-connectivity are different; a hypergraph with a single
edge containing every vertex is connected but not partition
connected.

%%%%
\subsection{Related work \labelsec{related}} Our results expand
theorems spanning graph theory, matroids and algorithms.  By
treating the problem in the most general setting, we will obtain
many of the results listed in this section as corollaries of our
more general results.

In this paragraph, we use {\em graph} in its usual sense, i.e.
as a $2$-uniform hypergraph.

\paragraph{Graph Theory and Rigidity Theory.}  Sparsity is
closely related to graph arborescence. The well-known results of
Tutte \cite{tutte61} and Nash-Williams \cite{Na61} show the
equivalence of
$(k,k)$-tight 
 graphs and graphs that can
be composed into $k$ edge-disjoint spanning trees.  A theorem of Tay
\cite{Ta1,Ta2} relates such graphs to generic rigidity of
bar-and-body structures in arbitrary geometric dimension. The
$(2,3)$-tight $2$-dimensional graphs play an important role in
rigidity theory. These are the generically minimally rigid graphs
\cite{laman70} (also known as Laman graphs), and have been studied
extensively. Results of Recski \cite{Re84I,Re84II} and Lovász
and Yemini \cite{lovasz:yemini} relate them to adding any edge to
obtain a $2$-arborescence. The most general results on $2$-graphs
were proven by Haas in \cite{haas:2002}, who shows the equivalence
of $(k,k+a)$-sparse graphs and graphs which decompose into $k$
edge-disjoint spanning trees after the addition of any $a$ edges. In
\cite{maps} Haas et al. extend this result to graphs that
decompose into edge-disjoint spanning maps, showing that
$(k,\ell)$-sparse graphs are those that admit such a map
decomposition after the addition of any $\ell$ edges.

For hypergraphs, Frank et al. study the $(k,k)$-sparse case in
\cite{frank2001}, generalizing the Tutte and Nash-Williams theorems
to partition connected hypergraphs.

\paragraph{Matroids.} Edmonds \cite{Ed65} used a matroid union approach to
characterize the $2$-graphs that can be decomposed into $k$ disjoint
spanning trees and described the first algorithm for recognizing
them.  White and Whiteley \cite{whiteley:matroids} first recognized
the matroidal properties of general $(k,\ell)$-sparse graphs.

In  \cite{whiteley:union-matroids}, Whiteley used a classical
theorem of Pym and Perfect \cite{pym-perfect} to show that the
$(k,\ell)$-tight $2$-graphs are exactly those that decompose into an
$\ell$-arborescence and $(k-\ell)$-map for $0\le \ell\le k$.

In the hypergraph setting, Lorea  \cite{lorea} described the first
generalization of graphic matroids to hypergraphs.  In
\cite{frank2001}, Frank et al. used a union matroid approach to
extend the Tutte and Nash-Williams theorems to arbitrary
hypergraphs.

\paragraph{Algorithms.} Our algorithms generalize the $(k,\ell)$-sparse
graph pebble games of Lee and Streinu \cite{LeSt05}, which in turn
generalize the pebble game of Jacobs and Hendrickson \cite{JaHe97}
for planar rigidity (which would be a $(2,3)$-pebble game in the
sense of \cite{LeSt05}). The elegant pebble game of \cite{JaHe97},
first analyzed for correctness in \cite{Be03}, was intended to be an
easily implementable alternative to the algorithms based on
bipartite matching discovered by Hendrickson in
\cite{hendrickson-thesis}.

The running time analysis of the $(2,3)$-pebble game in  \cite{Be03}
showed its running time to be dominated by $O(n^{2})$ queries about
whether two vertices are in the span of a rigid component.  This
leads to a data structure problem, considered explicitly in
\cite{LeSt05,cccg}, where it is shown that the running time of the
general $(k,\ell)$-pebble game algorithms on $2$-graphs is $O(n^2)$.

For certain special cases of $k$ and $\ell$, algorithms with better
running times have been discovered for $2$-multigraphs. Gabow and
Westermann \cite{GaWe88} used a matroid union approach to achieve a
running time of $O(n^{3/2})$ for the {\bf extraction} problem when
$\ell\le k$. They also find the set of edges that are in some
component, which they call the {\bf top clump}, with the same
running time as their extraction algorithm.  We observe that the
{\bf top clump} problem coincides with the components problem only
for the $\ell=0$ case. Gabow and Westermann also derive an
$O(n^{3/2})$ algorithm for the {\bf decision} problem for
$(2,3)$-sparse (Laman) graphs, which is of particular interest due
to the importance of Laman graphs in many rigidity applications.
Using a matroid intersection approach, Gabow \cite{gabow1995}
also gave an $O((m+n)\log n)$ algorithm for the extraction problem
for $(k,k)$-sparse $2$-graphs.

%%%%
\subsection{Our Results \labelsec{results}} We describe our results
in this section.

\paragraph{The structure of sparse hypergraphs.} We first describe
conditions for the existence of tight hypergraphs and analyze the
structure of the components of sparse ones. The theorems of this
section are generalizations of results from \cite{LeSt05,szego} to
hypergraphs of dimension $d\ge 3$.

\begin{theorem}[{\bf Existence of tight hypergraphs}]\labelthm{good-range}
There exists an $n_1$ depending on $s$, $k$ at $\ell$ such that
uniform tight graphs on $n$ vertices exist for all values of $n\ge
n_1$. In the smaller range $n<n_1$, such tight graphs may not exist.
\end{theorem}

\begin{theorem}[{\bf Block Intersection and Union}]
If $B_1$ and $B_2$ are blocks of a sparse graph $G$, $0\le \ell\le
ik$, and $B_1$ and $B_2$ intersect on at least $i$ vertices, then
$B_1\cup B_2$ is a block and the subgraph induced by $V(B_1)\cap
V(B_2)$ is a block. \labelthm{block-structure}
\end{theorem}

\begin{theorem}[{\bf Disjointness of Components}]
If  $C_1$ and $C_2$ are components of a sparse graph $G$, then
$E(C_1)$ and $E(C_2)$ are disjoint and $\card{V(C_1)\cap V(C_2)}<
s$. If $\ell\leq k$, then the components are vertex disjoint.  If
$\ell=0$, then there is only one component.
\labelthm{component-structure}
\end{theorem}

\paragraph{Hypergraph decompositions.} Extending the results of Tutte
\cite{tutte61}, Nash-Williams \cite{Na61}, Recski
\cite{Re84I,Re84II}, Lovász and Yemini \cite{lovasz:yemini},
 Haas et al. \cite{haas:2002,maps},
and Frank et al. \cite{frank2001}, we characterize the hypergraphs
that become $k$-arborescences after the addition of any $\ell$
edges.

\begin{theorem}[{\bf Generalized Lov{\'{a}}sz-Recski Property}]
Let $G$ be $(k,\ell)$-tight hypergraph with $\ell\ge k$.  Then the
graph $G'$ obtained by adding any $\ell-k$ edges of dimension at
least 2 to $G$ is a $k$-arborescence.
\labelthm{trees-after-adding-any}
\end{theorem}
In particular, the important special case in which $k=\ell$ was
proven by Frank et al. \cite{frank2001}.

\paragraph{Decompositions into maps.} We also extend the results of Haas
et al.  \cite{maps} to hypergraphs. This theorem can also be seen
as a generalization of the characterization of Laman graphs in
\cite{hendrickson-thesis}.

\begin{theorem}[{\bf Generalized Nash-Williams-Tutte Decompositions}]
A graph $G$ is a $k$-map if and only if $G$ is $(k,0)$-tight.
\labelthm{k-maps-are-tight}
\end{theorem}

\begin{theorem}[{\bf Generalized Haas-Lov{\'{a}}sz-Recski Property for Maps}]
The graph $G'$ obtained by adding any $\ell$ edges from $K_n^{k,0}-G$ to a
$(k,\ell)$-tight graph $G$ is a $k$-map.
\labelthm{maps-after-adding-any}
\end{theorem}

Using a matroid approach, we also generalize a theorem of Whiteley
\cite{whiteley:union-matroids} to hypergraphs.
\begin{theorem}[{\bf Maps and Trees Decomposition}]
Let $k\ge \ell$ and $G$ be tight.  Then $G$ is the union of an
$\ell$-arborescence and a $(k-\ell)$-map. \labelthm{maps-and-trees}
\end{theorem}

\paragraph{Pebble game constructible graphs.} The main theorem of this
paper, generalizing from $s=2$ in \cite{LeSt05} to hypergraphs of any
dimension, is that the matroidal families of sparse graphs coincide
with the pebble game graphs.

\begin{theorem}[{\bf Main Theorem: Pebble Game Constructible
Hypergraphs}]
 Let $k$, $\ell$, $n$ and $s$ meet the conditions of
\refthm{good-range}.  Then a hypergraph $G$ is sparse if and only if
it has a pebble game construction.
\labelthm{sparse-graphs-are-pebble-game-graphs}
\end{theorem}

\paragraph{Pebble game algorithms.}
We also generalize the pebble game {\em algorithms} of
\cite{LeSt05} to hypergraphs.  We present two algorithms, the {\bf
basic pebble game} and the {\bf pebble game with components}.

We show that on an  $s$-uniform input $G$ with $n$ vertices and $m$
edges, the basic pebble game solves the  {\bf decision} problem in
time $O((s+\ell)sn^2)$ and space $O(n)$.  The {\bf extraction}
problem is solved by the basic pebble game in time $O((s+\ell)dnm)$
and space $O(n+m)$.  For the {\bf optimization} problem, the basic
pebble game uses time $O((s+\ell)snm+m\log m)$ and space $O(n+m)$.

On an  $s$-uniform input $G$ with $n$ vertices and $m$ edges, the
pebble game with components solves  the {\bf decision}, {\bf
extraction}, and {\bf components} problems  in time
$O((s+\ell)sn^s+m)$ and space $O(n^s)$.  For the optimization
problem, the pebble game with components takes time
$O((s+\ell)sn^s+m\log m)$.

\paragraph{Critical representations.}
As an application of the pebble game,
we obtain lower-dimensional representations for certain classes of
sparse hypergraphs, generalizing a result from Lovász
 \cite{lovasz-representation} concerning lower-dimensional representations for (hypergraph) trees.

\begin{theorem}[{\bf Lower Dimensional and Critical Representations}]
\labelthm{representation} $G$ is a critical sparse hypergraph of
dimension $s$  if and only if the representation found by the pebble
game construction coincides with $G$. This implies that $G$ is
$s$-uniform and $\ell \leq sk-1$.
\end{theorem}

The proof of \refthm{representation} is based on a modified version of
the pebble game (described below) that solves the {\bf representation}
problem.  Its complexity is the same as that of the pebble game with
components: time $O((s+\ell)sn^s+m)$ and space $O(n^s)$ on an $s$-graph.

As corollaries to \refthm{representation}, we obtain:

\begin{corollary}[Lovász \cite{lovasz-representation}]
$G$ is an $s$-dimensional $k$-arborescence if and only if it is
represented by a 2-uniform $k$-arborescence $H$.
\end{corollary}

\begin{corollary}
$G$ is a $k$-map if and only if it is represented by a $k$-map with
edges of dimension  $1$.
\end{corollary}

\begin{corollary}
$G$ has a maps-and-trees decomposition if and only if $G$ is
represented by a graph with edges of dimension at most 2 that has a
maps-and-trees decomposition.
\end{corollary}

\section{The pebble game} The {\bf pebble game} is a family of
algorithms indexed by nonnegative integers $k$ and $\ell$.

The game is played by a single player on a fixed finite set of
vertices. The player makes a finite sequence of moves;  a move
consists of the addition and/or orientation of an edge. At any
moment of time, the state of the game is captured by a graph: we
call it a {\bf pebble game graph}.

Later in this paper, we will use the pebble game as the basis of
efficient algorithms for the computational problems defined above in
\refsec{sparse}.

We describe the pebble game in terms of its initial configuration
and the allowed moves.

{\bf Initialization:} in the beginning of the pebble game, $H$ has
$n$ vertices and no edges.  We start by placing $k$ pebbles on each
vertex of $H$.

{\bf Add edge:} Let $e\subset V$ be a set of vertices with at least
 $\ell+1$ pebbles on it.  Add $e$ to $E(H)$.  Pick up a pebble from
 any $v\in e$, and make $v$ the tail of $e$.

\begin{figure}%[htbp]
\centering
\subfigure[]{\includegraphics[width=2 in]{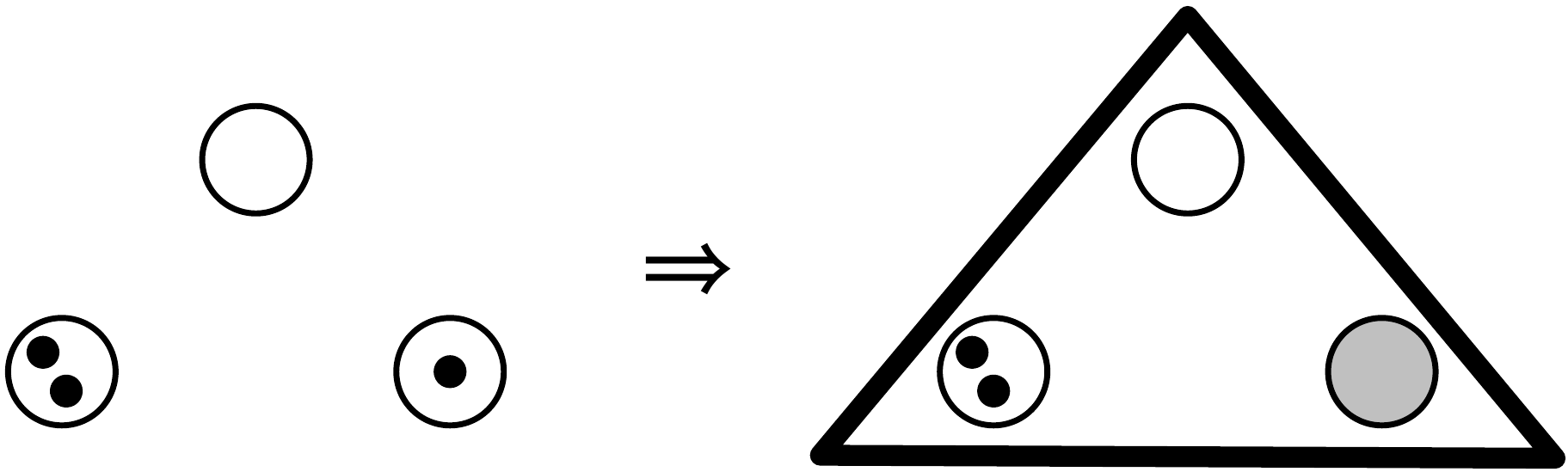}}
\hspace{.3 in}
\subfigure[]{\includegraphics[width=2 in]{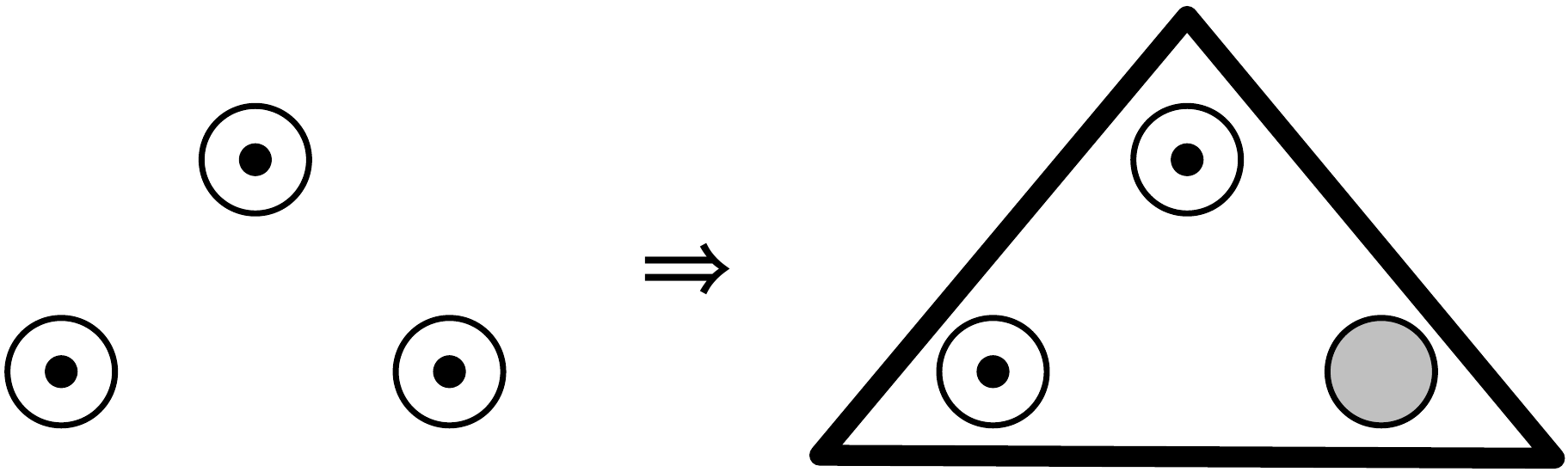}}
\subfigure[]{\includegraphics[width=2 in]{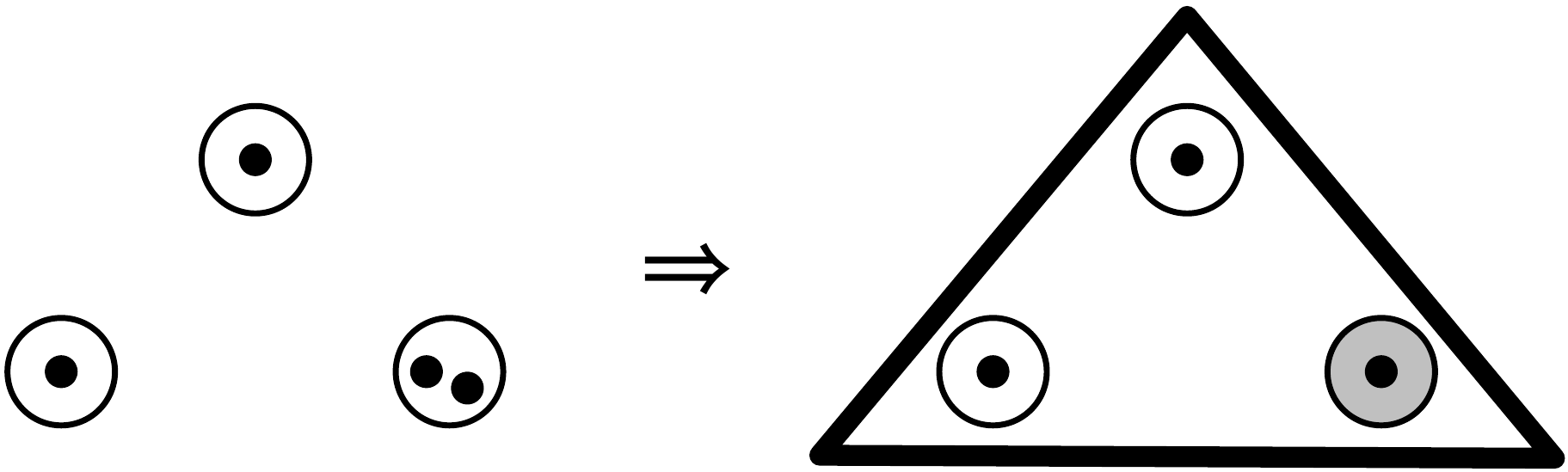}}
\caption{Adding a $3$-edge in the $(2,2)$-pebble game.  In all cases, the edge,
shown as a triangle, may be added because there are at least three pebbles present.
The tail of the new edge is  filled in; note that in (c) only one of the pebbles
on the tail is picked up.}
\label{fig.colored-add-edge}
\end{figure}

 \reffig{colored-add-edge} shows an example of this move in the $(2,2)$-pebble game.

 {\bf Pebble shift:} Let $v$ a vertex with at least one pebble on it,
 and let $e$ be an edge with $v$ as one of its ends, and with tail $w$.
 Move the pebble to $w$ and make $v$ the tail of $e$.

\begin{figure}%[htbp]
\centering
\subfigure[]{\includegraphics[width=2 in]{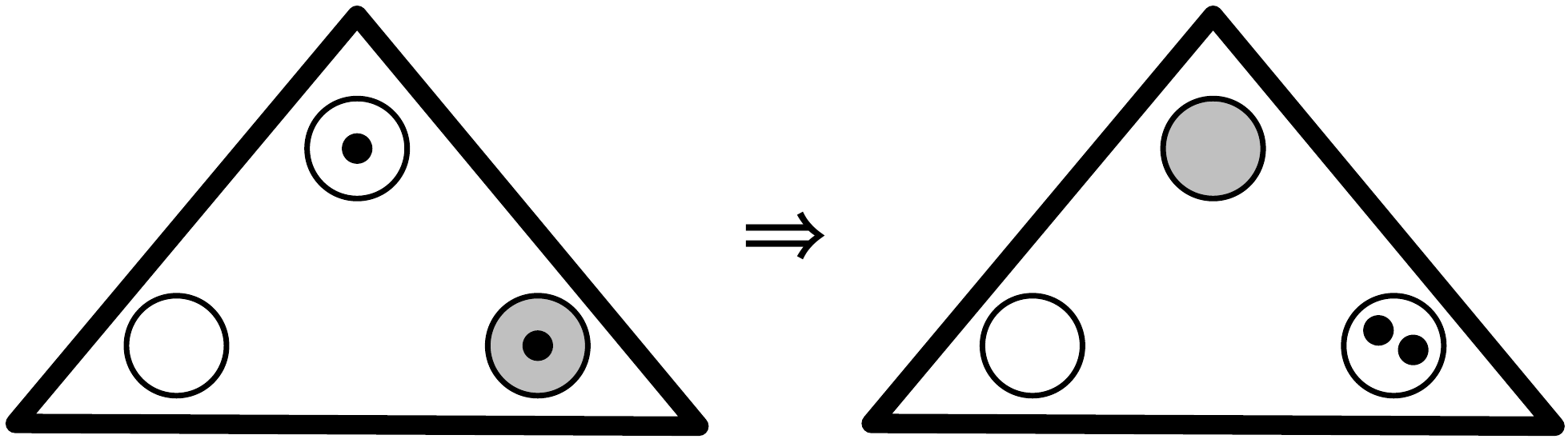}}
\hspace{.3 in}
\subfigure[]{\includegraphics[width=2 in]{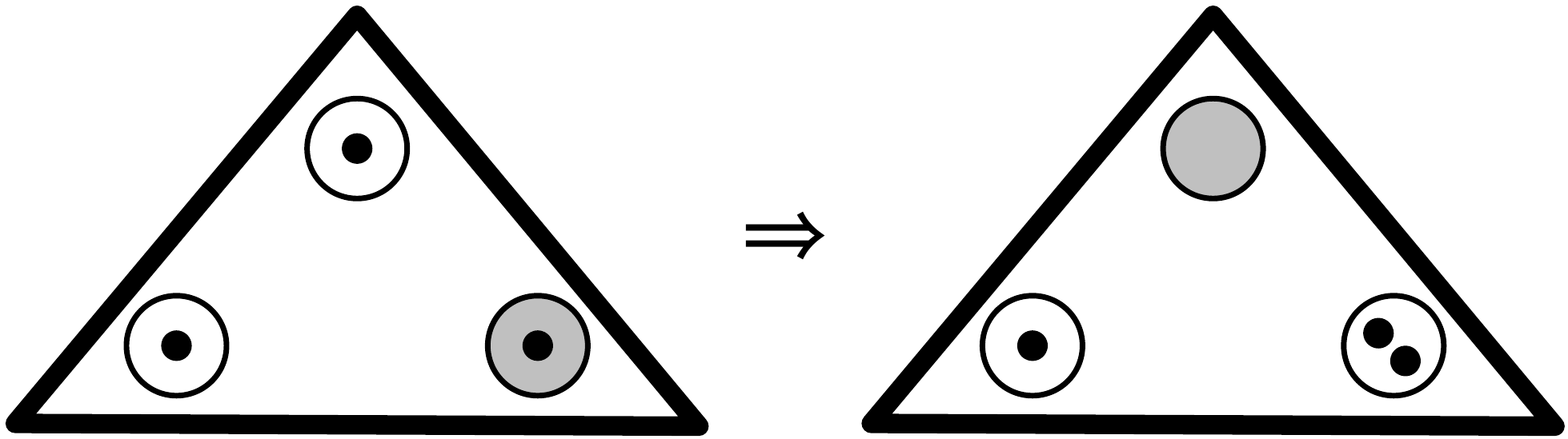}}
\caption{Moving a pebble along a $3$-edge in the $(2,2)$-pebble game.  The
tail of the edge is filled in.  Observe that in (b) the only change is to the
orientation of the edge and the location of the pebble that moved.}
\label{fig.colored-pebble-shift}
\end{figure}

 \reffig{colored-pebble-shift} shows an example of this move in the $(2,2)$-pebble game.

 The output of playing the pebble game is its complete configuration, which includes
 an oriented pebble game graph.

 {\bf Output:} At the end of the game, we obtain the oriented hypergraph $H$,
 and a map $\peb$ from $V$ to $\mathbb{N}$ such that for each vertex $v$,
 $\peb (v)$ is the number of pebbles on $v$.

\paragraph{Comparison to Lee and Streinu.}
The hypergraph pebble game extends the framework
developed in \cite{LeSt05} for $2$-graphs.
The main challenge was to come up with the concept of orientation of
hyperedges and of moving  the pebbles in a way that generalizes
depth-first search for $2$-graphs.  Specializing our
algorithm to $2$-uniform hypergraphs gives back the algorithm of
\cite{LeSt05}.

%%%%%%%
\section{Properties of sparse hypergraphs\labelsec{hypersparse}}

We next develop properties of sparse graphs, starting with the
conditions on $s$, $k$, $\ell$ and $n$ for which there are tight
graphs.

\begin{lemma}
If $\ell\ge ik$, and $G$ is sparse, then $s>i$.
\labellem{sparse-graph-rank}
\end{lemma}
\begin{proof}
If $i\ge s$, then for any edge $e$ of dimension $s$ the ends of $e$ are a set of
vertices for which \refeq{subset} fails.
\end{proof}

As an immediate corollary, we see that the class of uniform sparse graphs is
trivial when $\ell\ge sk$.
\begin{lemma}
If $\ell\ge sk$, then the class of $s$-uniform $(k,\ell)$-sparse
graphs contains only the empty graph.
\labellem{when-no-sparse-graphs}
\end{lemma}

We also observe that when $\ell<0$, the union of two disjoint sparse graphs need not
be sparse.  Since this is a desirable property, for the moment we focus on the case
in which $\ell\ge 0$.  Our next task is to further subdivide this range.

\begin{lemma}
Let $G$ be sparse and uniform.  The multiplicity of parallel edges
in $G$ is at most $sk-\ell$. \labellem{loops-and-parallel-edges}
\end{lemma}
\begin{proof}
\refeq{subset} holds for no more than $sk-\ell$ parallel edges of
dimension $s$.
\end{proof}

The next lemmas establish a range of parameters for which there are
tight graphs.
\begin{lemma}
Let $\ell\ge (s-1)k$.  There are no tight subgraphs on $n<s$ vertices.
\labellem{lower-trivial-range}
\end{lemma}
\begin{proof}
By \reflem{loops-and-parallel-edges} no sparse subgraph may contain edges
of dimension less than $s$.
\end{proof}

\begin{lemma}
If $\ell\ge (s-1)k$ then there is an $n_1$ depending on $s$, $k$ at
$\ell$ such that for $n\ge n_1$ there exist tight $s$-uniform graphs
on $n$ vertices. For $n<n_1$, there may not be tight uniform graphs.
\labellem{bad-range}
\end{lemma}
\begin{proof}
When $\ell\ge (s-1)k$ there are no loops in any sparse graph.  Also, by
\reflem{loops-and-parallel-edges} no edge in a uniform graph
has multiplicity greater
than $k$ in a sparse graph.  It follows that any tight uniform graph is a
subgraph of the complete $s$-uniform graph on $n$ vertices, allowing
edge multiplicity $k$.

For tight uniform subgraphs to exist, we need to have
\begin{eqnarray}
kn-\ell\le k\binom{n}{s}
\end{eqnarray}
Since the function $f(n)=kn^ss^{-s}-kn+\ell$ is asymptotically
positive, the desired $n_1$ must exist.

Notice that there is no tight $2$-uniform graph for $n=3$, $k=3$ and
$\ell=5$; the complete graph $K_3$ has only 3 edges, and by
\reflem{loops-and-parallel-edges} any $(3,5)$-sparse graph must be
simple. Such examples can be constructed for all values of $n\le
n_1$.
\end{proof}

We next turn to showing that tight graphs exist.
\begin{lemma}
Suppose that $\ell\ge (s-1)k$ and that $n\ge n_1$, where $n_1$ is taken
as in \reflem{bad-range}.  Then there are tight graphs on $n$ vertices.
\labellem{tight-graphs-exist}
\end{lemma}
\begin{proof}
Start with the complete $d$-uniform hypergraph
with $k$ parallel edges, $K^{k}_{n_1}$.  Identify a vertex $v$ and discard up to $\ell$
edges having $v$ as an end until the resulting graph $G_{n_1}$ is sparse.  This graph
must be sparse: any subgraph $H$ not spanning $v$ is sparse, as is any
subgraph containing only edges spanning $v$ by construction.
Since $G_{n_1}$ is maximally sparse, it is tight.

To complete the proof, proceed inductively: create $G_n$ from $G_{n-1}$ by
adding a new vertex and $k$ edges having the new vertex as an endpoint
such that the subgraph induced by the new edges is sparse.
\end{proof}

We next characterize the range of parameters for which there are
tight graphs.
\begin{restate}{good-range}{{\bf Existence of tight hypergraphs}}
    There is an $n_1$ depending on $s$, $k$ at $\ell$ such that
    for $n\ge n_1$ there are uniform tight graphs on $n$ vertices.
    For $n<n_1$, there may not be tight graphs.
\end{restate}
\begin{proof}%[Proof of \refthm{good-range}]
Immediate from
\reflem{bad-range} and \reflem{tight-graphs-exist}; the existence of tight uniform
hypergraphs implies the existence of tight hypergraphs.
\end{proof}

We next turn to the structure of blocks and components.
\begin{restate}{block-structure}{{\bf Block Intersection and Union}}
If $B_1$ and $B_2$ are blocks of a sparse graph $G$, $0\le \ell\le ik$,
$B_1$ and $B_2$ intersect on at least $i$
vertices, then $B_1\cup B_2$ is a block and the subgraph
induced by $V(B_1)\cap V(B_2)$ is a block.
%\labelthm{block-structure}
\end{restate}

\begin{proof}%[Proof of \refthm{block-structure}]
Let $m_i=\card{E(B_i)}$ for $i=1,2$; similarly let $v_i=\card{V(B_i)}$.
Also let $m_\cap=\card{E(B_1)\cap E(B_2)}$, $m_\cup=\card{E(B_1)\cup E(B_2)}$,
$v_\cup=\card{V(B_1)\cup V(B_2)}$, and
$v_\cap=\card{V(B_1)\cap V(B_2)}$.

The sequence of inequalities
\begin{eqnarray}
kn_\cup-\ell\ge m_\cup=m_1+m_2-m_\cap\ge kn_1-\ell+kn_2-\ell-kn_\cap+\ell=kn_\cup-\ell
\end{eqnarray}
holds whenever  $n_\cap\ge i$, which shows that $B_1\cup B_2$ is a block.

From the above, we get
\begin{eqnarray}
m_\cap=m_1+m_2-m_\cup=kn_1-\ell+kn_2-\ell-kn_\cup+\ell=
kn_\cap-\ell,
\end{eqnarray}
completing the proof.
\end{proof}

From \refthm{block-structure}, we obtain the first part of \refthm{component-structure}.
\begin{lemma}
If  $C_1$ and $C_2$ are components of a $(k,\ell)$-sparse graph $G$  then
$E(C_1)$ and $E(C_2)$ are disjoint and $\card{V(C_1)\cap V(C_2)}< s$.
\labellem{component-structure}
\end{lemma}
\begin{proof}
Observe that since $0\le \ell<sk$, components with non-empty edge intersection
are blocks meeting the condition of \refthm{block-structure}, as components
intersecting on $s$ vertices.  Since components are maximal, no two components
may meet the conditions of \refthm{block-structure}.
\end{proof}

For certain special cases, we can make stronger statements about the components.
\begin{lemma}
The components of a $(k,k)$-sparse graph are vertex disjoint.
\labellem{tree-components}
\end{lemma}
\begin{proof}
Observe that $\ell\le k$ and apply \refthm{block-structure} as above with $i=1$.
\end{proof}

\begin{lemma}
There is at most one component in a $(k,0)$-sparse graph.
\labellem{map-components}
\end{lemma}
\begin{proof}
Applying \refthm{block-structure} with $i=0$ shows that
the components of a
$(k,0)$-sparse graph are vertex disjoint.  Now suppose that
$C_1$ and $C_2$ are distinct components of a $(k,0)$-sparse
graph.  Then, using the notation of \refthm{block-structure},
$m_1+m_2=kn_1+kn_2=kn_\cup$, which implies that
$C_1\cup C_2$ is a larger component, contradicting the
maximality of $C_1$ and $C_2$.
\end{proof}

Together these lemmas prove the following result about the structure of components.
%\begin{theorem}
\begin{restate}{component-structure}{{\bf Disjointness of Components}}
If  $C_1$ and $C_2$ are components of a sparse graph $G$, then
$E(C_1)$ and $E(C_2)$ are disjoint and $\card{V(C_1)\cap V(C_2)}< s$.  If $k=\ell$, then the components are vertex disjoint.  If $\ell=0$, then there
is only one component.
%\labelthm{component-structure}
\end{restate}
\begin{proof}%[Proof of \refthm{component-structure}]
Immediate from \reflem{component-structure}, \reflem{tree-components},
and \reflem{map-components}.
\end{proof}

\section{Hypergraph Decompositions}

In this section we investigate links between tight hypergraphs and
decompositions into edge-disjoint maps and trees.

\subsection{Hypergraph arboricity} We now generalize
results of Haas \cite{haas:2002} and Frank et al. \cite{frank2001}
to prove an equivalence between
sparse hypergraph and those for which adding any $a$ edges results
in a $k$-arborescence.

We will make use of the following important result from  \cite{frank2001}.
\begin{proposition}[ Frank et al.  \cite{frank2001} ]
A hypergraph $G$ is a $k$-arborescence if and only if $G$ is $(k,k)$-tight.
\labelprop{frank-k-arborescence}
\end{proposition}

%\begin{theorem}
\begin{restate}{trees-after-adding-any}{{\bf Generalized Lovász-Recski Property}}
Let $\ell\ge k$ and let $G$ be tight.  Then the graph $G'$ obtained by adding any
$\ell-k$ edges of dimension at least 2 to $G$ is a $k$-arborescence.
%\labelthm{trees-after-adding-any}
%\end{theorem}
\end{restate}
\begin{proof}%[Proof of \refthm{trees-after-adding-any}]
Suppose that $G$ is tight and that $\ell\ge k$.  Let $G'=(V,F)$ be a graph obtained by adding
$\ell-k$ edges of dimension at least 2 to $G$, and consider a subset $V'$ of $V$.  It
follows that
\begin{eqnarray}
\card{E_{G'}(V')}\le \card{V'}+\ell-k\le kn'-\ell+\ell-k=kn'-k,
\end{eqnarray}
which implies that $G'$ is $(k,k)$-tight,  since $\card{F}=kn-k$. By
\refprop{frank-k-arborescence} $G'$ is a $k$-arborescence.

Conversely, if adding any $\ell-k$ edges to $G$ results in a $(k,k)$-tight
graph, then $G$ must be tight; if $V'$ spans more than $kn-\ell$ edges in $G$, then
adding $\ell-k$ edges to the the span of $V'$ results in a graph which is not
$(k,k)$-sparse.
\end{proof}

\subsection{Decompositions into maps} The main result of  this
section shows the equivalence of the $(k,0)$-tight graphs and
$k$-maps.  As an application, we obtain a characterization of all
the sparse hypergraphs in terms of adding {\em any} edges.

\begin{restate}{k-maps-are-tight}{{\bf Generalized Nash-Williams-Tutte Decompositions}}
A graph $G$ is a $k$-map if and only if $G$ is $(k,0)$-tight.
%\labelthm{k-maps-are-tight}
\end{restate}
\begin{proof}%[Proof of \refthm{k-maps-are-tight}]
Let $G=(V,E)$ be a hypergraph with $n$ vertices and $kn$ edges.
Let $B^k_G=(V_k,E,F)$ be the bipartite graph with one vertex class indexed by $E$ and
the other by $k$ copies of $V$.  The edges of $B^k_G$ capture the incidence structure of
$G$.  That is, we define $F=\{ v_ie : e=vw, e\in E, i=1,2,\ldots,k\}$; i.e., each edge vertex in $B$
is connected to the $k$ copies of its endpoints in $B^k_G$.  \reffig{k3-bipartite-example} shows $K_3$
and $B^1_{K_3}$.

\begin{figure}[htbp]
\centering %%
\includegraphics[height=1.0 in]{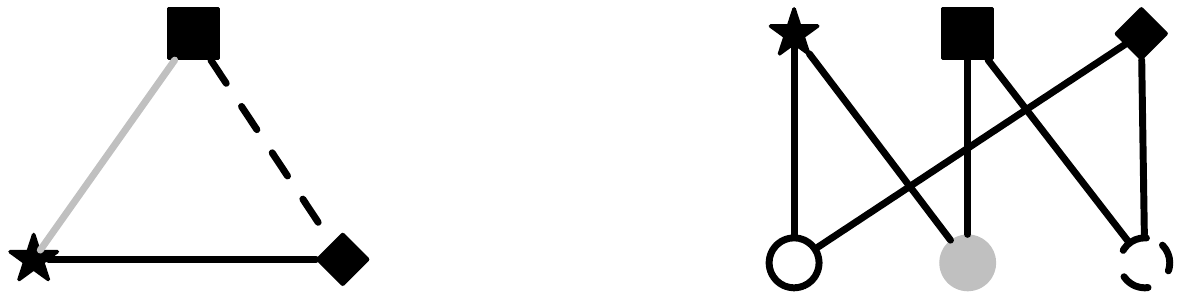}
\caption{The $(1,0)$-sparse 2-graph $K_3$ and its associated bipartite graphs
$B^1_{K_3}$.  The vertices and edges of $K_3$ are
matched to the corresponding vertices in $B^1_{K_3}$ by shape and
line style.}
\label{fig.k3-bipartite-example}
\end{figure}

Observe that
for any subset $E'$ of $E$,
\begin{eqnarray}
\card{N_{B^k_G}(E)}=k\card{V(E')}\ge \card{E}.
\labeleq{map-hall-condition}
\end{eqnarray}
if and only if $G$ is $(k,0)$-sparse.  By Hall's theorem, this implies that
$G$ is $(k,0)$-tight if and only if $B^k_G$ contains a perfect matching.

\fig{width=4 in}{An orientation of a 2-dimensional $2$-map $G$ and the associated bipartite matching in $B^2_G$.}{bipartite-2-map}
The edges matched to the $i$th copy of $V$ correspond to
the $i$th map in the $k$-map, as
shown for a 2-map in \reffig{bipartite-2-map}.
Assign as the tail of  each edge away from the vertex to which it is matched.
It follows that
each vertex has out degree one in the spanning subgraph matched to
each copy of $V$ as desired.
\end{proof}

\refthm{k-maps-are-tight} implies \refthm{maps-after-adding-any}.

%\begin{theorem}
\begin{restate}{maps-after-adding-any}{{\bf Generalized Haas-Lov{\'{a}}sz-Recski Property for Maps}}
    The graph $G'$ obtained by adding any $\ell$ edges from $K_n^{k,0}-G$ to a
    $(k,\ell)$-tight graph $G$ is a $k$-map.
%\labelthm{maps-after-adding-any}
\end{restate}
\begin{proof}%[Proof of \refthm{maps-after-adding-any}]
Similar to the proof of \refthm{trees-after-adding-any}.  Because the added edges come from
$K_n^{k,0}-G$, the resulting graph must be sparse.
\end{proof}

We see from the proof of \refthm{maps-after-adding-any}, that the condition of adding edges of dimension
at least 2 in \refthm{trees-after-adding-any} is equivalent to saying that the added edges
come from $K_n^{k,k}$.

To prove \refthm{maps-and-trees}, we need several results from matroid
theory.

\begin{proposition}
    Let $r$ be a non-negative, increasing, submodular set function on a finite
    set $E$.  Then the class
    $\mathcal{N}=\{A\subset E : \card{A'}\le r(A'), \forall A'\subset A \}$
    gives the independent sets of a matroid.
\end{proposition}
We say that $\mathcal{N}$ is generated by $r$.  In particular, we see
that our matroids of sparse hypergraphs are generated
by the function $r_{k,\ell}(E')=k\card{V(E')}-\ell$.

Pym and Perfect \cite{pym-perfect}
proved the following result about unions of such matroids.
\begin{proposition}[Pym and Perfect \cite{pym-perfect}]
Let $r_1$ and $r_2$ be non-negative, submodular, integer-valued
functions, and let $\mathcal{N}_1$ and $\mathcal{N}_2$
be matroids they generate.  Then the matroid union of
$\mathcal{N}_1$ and $\mathcal{N}_2$ is generated by $r_1+r_2$. \labelprop{pym}
\end{proposition}

Let $\mathcal{M}_{1,0}$ and $\mathcal{M}_{1,1}$ be the matroids which
have as bases the $(1,0)$-tight and $(1,1)$-tight hypergraphs
respectively.  That these are matroids is a result of White and Whiteley from
\cite{whiteley:matroids} proven in the appendix of this paper for completeness.
\refthm{k-maps-are-tight} and \refprop{frank-k-arborescence} imply that the bases
of these matroids are the maps and trees and that these matroids are
generated by the functions $r_{1,0}(E')=\card{V(E')}$ and
$r_{1,1}(E')=\card{V(E')}-1$.

With these observations we can prove \refthm{maps-and-trees}.

\begin{restate}{maps-and-trees}{{\bf Decompositions into maps and trees}}
Let $k\ge \ell$ and $G$ be tight.  Then $G$ is the union of an
$\ell$-arborescence and a $(k-\ell)$-map.
%\labelthm{maps-and-trees}
\end{restate}
\begin{proof}%[Proof of \refthm{maps-and-trees}]
We first observe that $r_{1,0}$
meets the conditions of \refprop{pym}.  Since $r_{1,1}$ does not (it
is not non-negative), we switch to the submodular function
\begin{eqnarray}
r'(V')=n'-c
\end{eqnarray}
where $c$ is the number of non-trivial partition-connected components
spanned by $V'$.  It follows that $r'$ is non-negative, since a
graph with no edges has no non-trivial partition-connected
components. Observe also, that if $V'$ spans $c$ partition-connected
components with $n_1,n_2,\ldots,n_c$ vertices we have
\begin{eqnarray}
r_{1,1}(V')=\sum_{i=1}^c(n_i-1)=n'-c=r'(V'),
\end{eqnarray}
since the partition-connected components are blocks of trees, and
thus disjoint.

Applying \refprop{pym} to $r_{1,0}$ and $r'$ now shows that
the union matroid of $k-\ell$ maps and $\ell$ trees is generated by
\begin{eqnarray}
r(V')=(k-\ell)r_{1,0}(V')+\ell r'(V')=(k-\ell)n'+\ell
n'-\ell,
\end{eqnarray}
proving that the union of the matroid with bases that decompose into
$(k-\ell)$ maps and $\ell$ trees is $\mathcal{M}_{k,\ell}$ as
desired.
\end{proof}

\section{Pebble game constructible graphs}
The main result of this section is that the matroidal sparse graphs are exactly
the ones that can be constructed by the pebble game.

We begin by
establishing some invariants that hold during the execution of the
pebble game.

\begin{lemma}During the execution of the pebble game, the following
invariants are maintained in $H$:
\begin{description}
\item [{\bf (I1)}] There are at least $\ell$ pebbles on $V$.
\item [{\bf (I2)}] For each vertex $v$, $\grsp v + \out v + \peb v=k$.
\item [{\bf (I3)}] For each $V'\subset V$, $\grsp V'+\out V'+\peb V'=kn'$.
\end{description}
\labellem{pebble-game-invariants}
\end{lemma}
\begin{proof}
{\bf (I1)} The number of pebbles on $V$ changes only after an {\bf add edge move}.
When there are fewer than $\ell+1$ pebbles, no {\bf add edge} moves are possible.

{\bf (I2)} This invariant clearly holds at the initialization of the pebble game.
We verify that each of the moves preserves {\bf (I2)}.  An {\bf add edge} move
consumes a pebble from exactly one vertex and adds one to its out degree or span.
Similarly, a {\bf pebble shift} move adds one to the out degree of the source and
removes a pebble while adding one pebble to the destination and decreasing its
out degree by one.

{\bf (I3)}  Let $V'\subset V$ have $n'$ vertices and span $m^{+}$ edges with
at least two ends.  Then
\begin{eqnarray}
\out V'=\sum_{v\in V'} \out v-m^{+}
\end{eqnarray}
and
\begin{eqnarray}
\grsp V'=m^{+}+\sum_{v\in V'}\grsp v.
\end{eqnarray}
Then we have
\[
\begin{split}
\grsp V'+ \out V'+\peb V' \\
&=\sum_{v\in V'} \out v-m^{+} + m^{+} +
\sum_{v\in V'} \grsp v + \sum_{v\in V'} \peb v \\
&=\sum_{v\in V'}(\out v+\grsp v+\peb v)  = kn',
\end{split}
\]
where the last step follows from {\bf (I2)}.
\end{proof}

From these invariants, we can show that the pebble game constructible graphs
are sparse.

\begin{lemma}Let $H$ be a hypergraph constructed with the pebble game.  Then $H$
is sparse.  If there are exactly $\ell$ pebbles on $V(H)$, then $H$ is tight.
\labellem{pebble-graphs-are-sparse}
\end{lemma}
\begin{proof}
Let $V'\subset V$ have $n'$ vertices and consider the configuration of the
pebble game immediately after the most recent {\bf add edge} move that
added to the span of $V'$.  At this point, $\peb V'\ge \ell$.  By \reflem{pebble-game-invariants}
{\bf (I3)},
\begin{eqnarray}
kn'\ge \grsp V'+\out V'+\ell.
\end{eqnarray}
When $\grsp V'>kn'-\ell$, this implies that $-1\ge\out V'$, which is a contradiction.

In the case where there are exactly $\ell$ pebbles on $V(H)$, \reflem{pebble-game-invariants}
{\bf (I3)} implies that $\grsp V=kn-\ell$.
\end{proof}

We now consider the reverse direction: that all the sparse graphs
admit a pebble game construction.  We start with the observation
that if there is a path in $H$ from $u$ to $v$, then if $v$ has a pebble
on it, a sequence of {\bf pebble shift} moves can bring the
pebble to $u$ from $v$.

\newcommand{\reach}{\ensuremath{\operatorname{reach}}}
Define the {\bf reachability region} of a vertex $v$ in $H$ as the set
\begin{eqnarray}
\reach{v}=\{u\in V : \text{there is a path in $H$ from $v$ to $u$}\}.
\end{eqnarray}

\begin{lemma}
Let $e$ be a set of vertices such that $H+e$ is sparse.  If $\peb e<\ell+1$,
then a pebble not on $e$ can be brought to an end of $e$.
\labellem{can-bring-another-pebble}
\end{lemma}
\begin{proof}
Let $V'$ be the union of the reachability regions of the ends of $e$; i.e.,
\begin{eqnarray}
V'=\bigcup_{v\in e}\reach v.
\end{eqnarray}
Since $V'$ is a union of reachability regions, $\out V'=0$.  As $H+e$
is sparse and $e$ is in the span of $V'$, $\grsp V'<kn'-\ell$.

It follows by \reflem{pebble-game-invariants} {\bf (I3)}, that $\peb V'\ge \ell+1$,
so there is a pebble on $V'-e$.  By construction there is a $v\in e$ such that
the pebble is on a vertex $u\in \reach v-e$.  Moving the pebble from $u$
to $v$ does not affect any of the other pebbles already on $e$.
\end{proof}

It now follows that any sparse hypergraph has a pebble game construction.
\begin{restate}{sparse-graphs-are-pebble-game-graphs}{{\bf The Main Theorem: Pebble Game Constructible Hypergraphs}}
Let $G$ be a $(k,\ell)$-sparse hypergraph with $k$, $\ell$ and $s$
meeting the conditions of \refthm{good-range}.  Then $G$ can be constructed by the pebble
game.
\end{restate}
\begin{proof}
For each edge $e$ of  $G$ in any order, inductively apply \reflem{can-bring-another-pebble} to
the ends of $e$ until there are $\ell+1$ of them.  At this point, use an {\bf add edge} move to
add $e$ to $H$.
\end{proof}

It is instructive to note that the pebble game invariants enforce
the matroid properties of the sparse graphs.  The $\ell+1$
acceptance condition enforces the constraints on $k$, $\ell$ and
$s$, and the proof of \reflem{can-bring-another-pebble} shows that
the order in which edges of a sparse graph are added does not matter
in a pebble game construction.

\section{Pebble games for Components and Extraction}
Until now we were concerned with characterizing sparse and tight
graphs.  In this section we describe efficient algorithms based on
pebble game constructions.

\subsection{The basic pebble game}
In this section we develop the basic $(k,\ell)$-pebble game for hypergraphs to solve the
{\bf decision} and {\bf extraction} problems.  We first describe the algorithm.

\begin{algorithm}[The $(k,\ell)$-pebble game]
$\quad$ \\
{\bf Input:} A hypergraph $G=(V,E)$ \\
{\bf Output:} `sparse', `tight' or `dependent.' \\
{\bf Method:} Initialize a pebble game construction on $n$
vertices.

For each edge $e$, try to collect $\ell+1$ pebbles on the ends of $e$.
Pebbles can be collected using depth-first search to find a path
to a pebble and then a sequence of {\bf pebble shift} moves to move it.

If it is possible to collect $\ell+1$ pebbles,
use an {\bf add edge} move to add $e$ to $H$.

If any edge was not added to $H$, output `dependent'.
If every edge was added and there are exactly
 $\ell$ pebbles left, then output `tight'.  Otherwise output `sparse'.
\labelalg{basic-pebble-game}
\end{algorithm}

\begin{figure}[htbp]
\centering %%
\subfigure[]{\label{fig.pebble-1}\includegraphics[height=1.5 in]{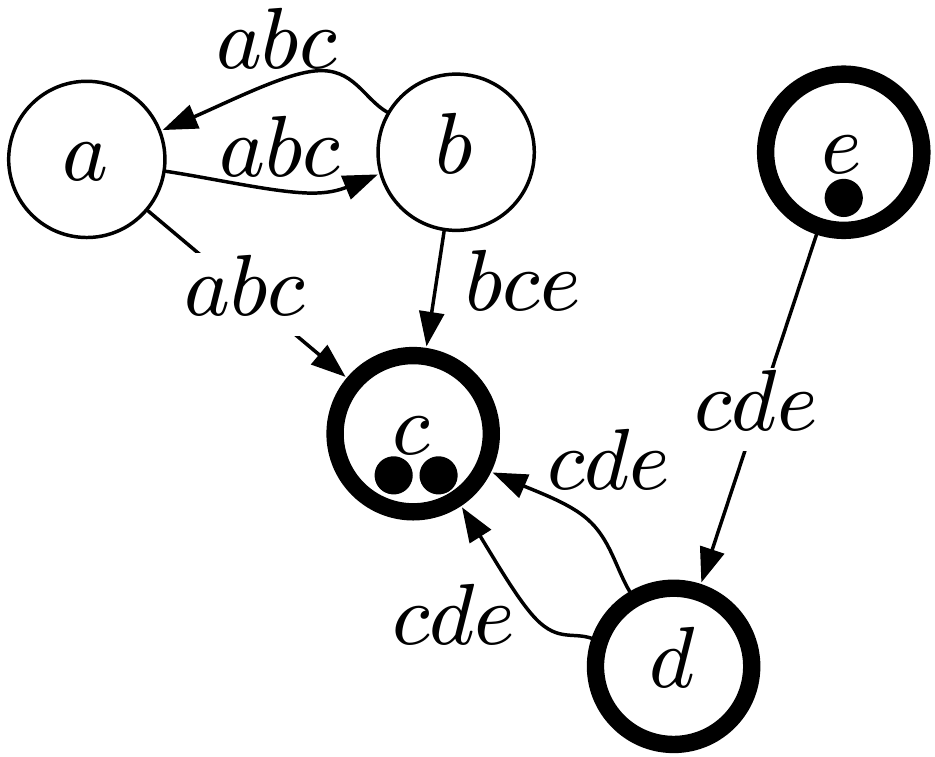}}
%\hspace{.3in}
\subfigure[]{\label{fig. pebble-2}\includegraphics[height=1.5 in]{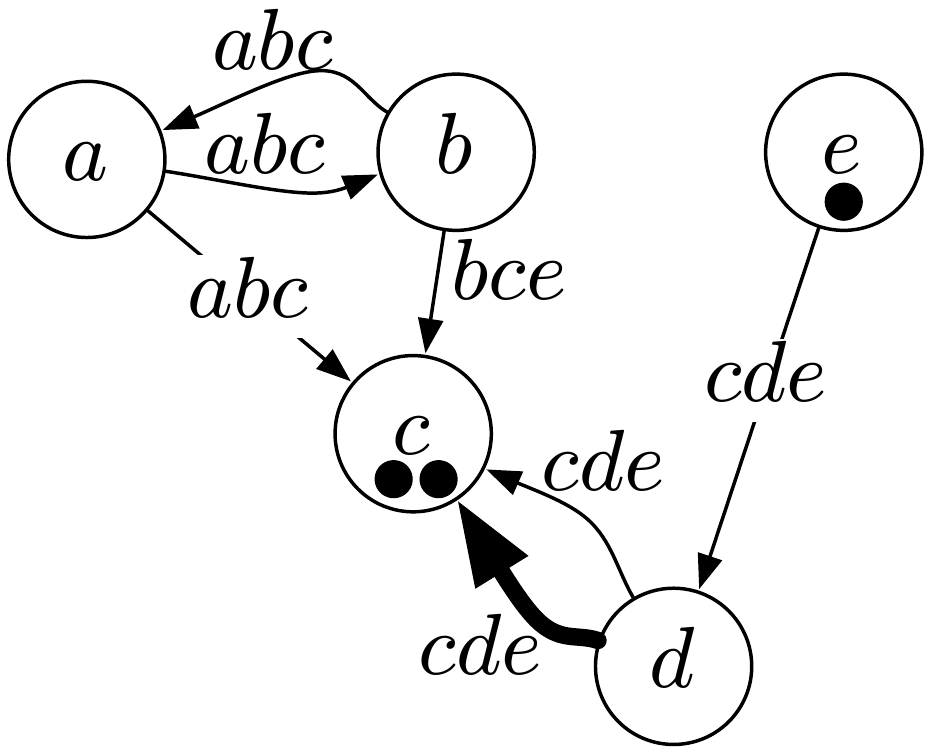}}
%\hspace{.3 in}
\subfigure[]{\label{fig. pebble-4}\includegraphics[height=1.5 in]{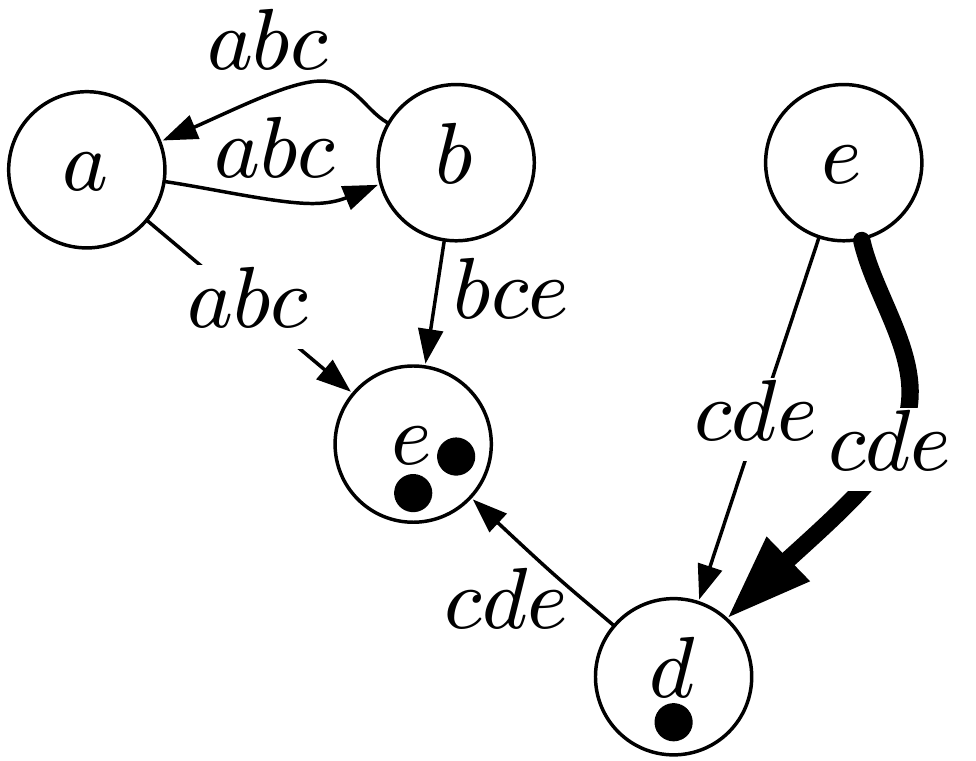}}
%\hspace{.3in}
%%
\subfigure[]{\label{fig. pebble-7}\includegraphics[height=1.5 in]{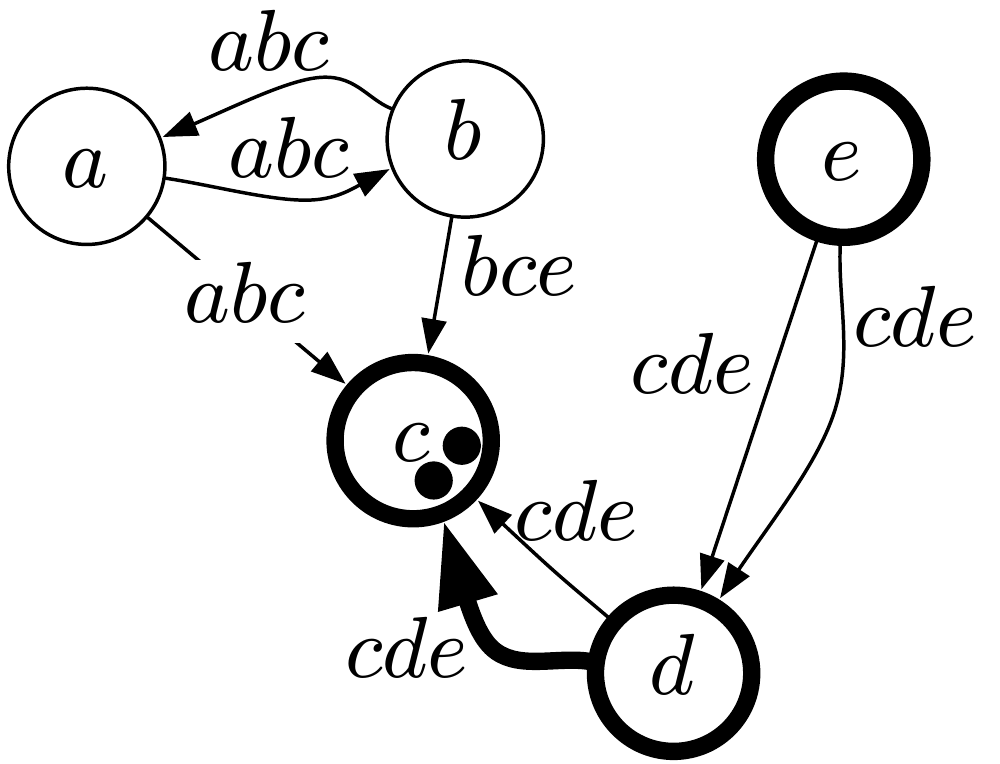}}
\caption{Collecting a pebble and accepting an edge in a $(2,2)$-pebble game on a
$3$-uniform hypergraph $H$.  $H$ is shown via a 2-uniform representation.
In (a), the edge being tested, $cde$ is shown with thick circles around the vertices.
 The pebble game starts a search
to bring a pebble to $d$.  (This choice is arbitrary; had $e$ been chosen first, the
edge would be immediately accepted.)
In (b) a path from $d$ to $e$ across the edge marked with a think line is found.
In (c)  the pebble is moved and the path is reversed; the new tail
of the edge marked with a think line is $e$.
In (d) the pebble is picked up, and the edge being checked is accepted.  The
tail of the new edge, marked with a thick line, is $d$.
}
\label{fig.collect-pebble}
\end{figure}

\reffig{collect-pebble} shows an example of collecting a pebble and accepting an edge.

The correctness of the basic pebble game for the {\bf decision}
and {\bf extraction} problems
 follows immediately from \refthm{sparse-graphs-are-pebble-game-graphs}.
For the {\bf optimization }
problem, sort the edges in order of increasing weight before starting;
the correctness follows from
\refthm{kl-matroid} and the  characterization of
matroids by the greedy algorithm (discussed in, e.g.,
\cite{oxley:matroid}).

The running time of the pebble game is dominated by the time needed
to collect pebbles.  If the maximum edge size in the hypergraph is
$s^*$, the time for one depth-first search is $O(s^*n+m)$, from
which it follows that the time to find one pebble in $H$ is
$O(s^*n)$. To check an edge requires no more than $s^*+\ell+1$
pebble searches, and $m$ edges need to be checked.  To summarize, we
have proven the following.

\begin{lemma}
Let $G$ be a hypergraph with $n$ vertices, $m$
edges, and maximum edge size $s^*$. The running time of the basic pebble game is
$O((s^*+\ell)s^*nm)$; for the {\bf decision} problem, this is
$O((s^*+\ell)s^*n^2)$, since $m=O(n)$.
\labellem{pebble-game-running-time}
\end{lemma}

All of the searching, marking, and pebble counting
can be done with $O(1)$ space per vertex.
Since  $H$ has $O(n)$ edges, the space complexity of the basic pebble game
is dominated by the size of the input.

\begin{lemma}
The space complexity of the basic pebble game is $O(m+n)$, where $m$ and $n$ are,
 respectively,
 the number of edges and vertices in the input.
 \labellem{pebble-game-space}
\end{lemma}

Together the preceding lemmas complete the complexity analysis.  The running time for the {\bf decision}
problem on a $d$-uniform hypergraph with $n$ vertices and $kn-\ell$ edges
is $O((s+\ell)sn^2)$,  and the space used $O(n)$.  For the {\bf optimization } problem, the running
time increases to $O((s+\ell)sn^2+n\log n)$ because of the sorting phase.

The {\bf extraction} problem is solved in time $O((s+\ell)snm)$ and space $O(n+m)$.

% YYYYYYY

\subsection{Detecting components}
In the next several sections we extend the basic pebble game to
solve the {\bf components} problem. Along the way, we also improve
the running time for the {\bf extraction} problem by developing a
more efficient way of discarding dependent edges. As the proof of
\reflem{pebble-game-running-time} shows, the time spent trying to
bring pebbles to the ends of dependent edges can be $\Omega(n^2)$ if
the edges are very large. We will reduce this to $O(s)$, improving
the running time.

We first present an algorithm to detect components.

\begin{algorithm}[Component detection]
    $\quad$ \\
{\bf Input:} An oriented hypergraph $H$ and $e$, the most recently accepted edge. \\
{\bf Output:} The component spanning $e$ or `free.' \\
{\bf Method:} When the algorithm starts, there are $\ell$ pebbles on the ends of $e$,
and a vertex $w$ is the tail of $e$.  If there are any other pebbles on $\reach w$,
stop and output `free.'  Otherwise let $C=\reach w$, and enqueue any vertex that is
an end of an edge pointing into $C$.

While there are more vertices in the queue, dequeue a vertex $u$.  If the only pebbles in
$\reach u$ are the $\ell$ on $e$, add $\reach w$ to $C$ and
enqueue any newly discovered vertex that is an end of an
edge pointing into $C$.

Finally, output $C$.
\labelalg{detect-components}
\end{algorithm}

In the rest of this section we analyze the correctness and running time of
\refalg{detect-components}.  We put off a discussion of the space required to
maintain the components until the next section.

We start with a technical lemma about blocks.

\begin{lemma}
Let $G$ be tight and $\ell>0$.  Then $G$ is connected.
\labellem{when-blocks-are-connected}
\end{lemma}
\begin{proof}
Consider a partition of $V$ into two subsets.  These span at most
$kn-2\ell$ edges by sparsity, but $G$ has $kn-\ell$ edges.
\end{proof}

\begin{lemma}
If \refalg{detect-components} outputs `free,' then $e$ is not spanned by any component.  Otherwise
the output $C$ of \refalg{detect-components} is the component spanning $e$.
\labellem{detect-components-correctness}
\end{lemma}
\begin{proof}
\refalg{detect-components} outputs `free' only when it is possible to collect at least $\ell+1$ pebbles
on the ends of $e$.  \reflem{pebble-graphs-are-sparse} shows that in this case, $e$ is not spanned by
any block in $H$ and thus no component.

Now suppose that \refalg{detect-components} outputs a set of vertices $C$.  By construction,
the number of free pebbles on $C$ is $\ell$.  Also, since $C$ is the union of reachability
regions, it has no out edges.  By \reflem{pebble-graphs-are-sparse}, $C$ spans a block in $H$.
Since \refalg{detect-components} does a breadth first search in $H$, $C$ is a maximal
connected block.

There are now two cases to consider.  When $\ell>0$, blocks are  connected
by \reflem{when-blocks-are-connected}.  If $\ell=0$, blocks may not be connected, but there
is only one component in $H$ by \reflem{map-components}; add $C$ to the component being maintained.
\end{proof}

For the running time of \refalg{detect-components} we observe that
$O(s^*)$ time is spent processing the vertices of each edge pointing
into $C$ for enqueueing and dequeuing.   Vertices are explored by
pebble searches only once; mark vertices accepted into $C$ and also
those from which pebbles can be reached to cut off the searches.
Since $H$ is $(k,\ell)$-sparse, it has $O(n)$ edges.  Summarizing,
we have shown the following.

\begin{lemma}
The running time of \refalg{detect-components} is $O(s^*n)$.
\labellem{detect-components-running-time}
\end{lemma}

%%%%
\subsection{The pebble game with components}
We now present an extension of the basic pebble game that solves the components
problem.

\begin{algorithm}[The $(k,\ell)$-pebble game with components]
    $\quad$ \\
{\bf Input:} A hypergraph $G=(V,E)$ \\
{\bf Output:} `Strict', `tight' or `dependent.' \\
{\bf Method:} Modify \refalg{basic-pebble-game} as follows.  When processing an
edge $e$ first check if it is spanned by a component.  If it is, then reject it.  Otherwise
collect $\ell+1$ pebbles on $e$ and accept it.  After accepting $e$, run
\refalg{detect-components} to find a new component if once has
been created.

Output the components discovered along with the output of the basic pebble
game.
\labelalg{components-pebble-game}
\end{algorithm}

The correctness of \refalg{components-pebble-game} follows
 from the fact that $H+e$ is sparse if and only if $e$ is not in the
 span of any component and \refthm{sparse-graphs-are-pebble-game-graphs}.

\begin{lemma}
\refalg{components-pebble-game} solves the {\bf decision}, {\bf extraction}
and {\bf components} problems.
\labelthm{components-pebble-game-is-correct}
\end{lemma}

%%%%
\subsection{Complexity of the pebble game with components}

We analyze the running time of the pebble game with components
in two parts: component maintenance and edge processing.

For component maintenance, we easily generalize the union pair-find data structures described in
 \cite{cccg}.  If $s^*$ is the largest size of an edge in $G$,
the complexity of
checking whether an edge is spanned by a component is $O(s^*)$, and
the total time spent updating the components discovered is $O(n^{s^*})$.
The complexity is dominated by maintaining a table with $n^{s^*}$ entries that
records with $s^*$-tuples are spanned by some component.

The time spent processing dependent edges is $O(s^*n^{s^*})$; they are exactly those
edges spanned by a component.  For each accepted edge, we need to
collect $\ell+1$ pebbles. The analysis is similar to that for the basic pebble game.
Since there are $O(n)$ edges accepted, we have the following total running time.

\begin{lemma}
The running time of \refalg{components-pebble-game} on a $s$-dimensional
hypergraph with $n$ vertices and $m$ edges is $O((s^*+\ell)s^*n^{s^*}+m)$.
\labellem{components-running-time}
\end{lemma}

Since the data structure used to maintain the components uses a table of
size $\Theta(n^{s^*})$, the space complexity of the pebble game with
components is the same on any input.

\begin{lemma}
The pebble game with components uses $O(n^s)$ space.
\labelthm{components-space}
\end{lemma}

Together the preceding lemmas complete the complexity analysis of the
pebble game with components.  The running time on an $s$-graph with
$n$ vertices and $m$ edges is $O((s+\ell)sn^s+m)$ and the space used is
$O(n^s)$.  For the optimization problem, the sorting phase of the
greedy algorithm takes an additional $O(m\log m)$ time.

% ZZZZZZZZ

\section{Critical representations}
As an application of the pebble game, we investigate the
circumstances under which we may represent a sparse hypergraph with
a lower dimensional sparse hypergraph.  The main result of this
section is a complete characterization of the critical sparse
hypergraphs for any $k$ and $\ell$.

Clearly, by \reflem{sparse-graph-rank}, when $\ell\ge (s-1)k$, every sparse
$s$-uniform hypergraph must be critical.  In this section we show that
these are the only $s$-uniform critical sparse hypergraph and describe an
algorithm for finding them.

We first present a modification of the pebble game to compute a
representation.  Only the {\bf add edge} and {\bf pebble shift}
moves need to change.

{\bf Represented add edge:} When adding an edge $e$ to $H$, create a set $r(e)$
which is the set of vertices with the $\ell+1$ pebbles used to certify that
$e$ was independent.

{\bf Represented pebble shift:} When a {\bf pebble shift} move makes an end
$v\notin r(e)$ the tail of $e$, add $v$ to $r(e)$ and remove any other element of $r(e)$.

Let $R$ be the oriented hypergraph with the edge set $r(e)$ for $e\in E(H)$.

We now consider the invariants of the represented pebble game.

\begin{lemma}\labellem{represented-pebble-game-invariants}
    The invariants {\bf (I1)}, {\bf (I2)}, and {\bf (I3)} hold in $R$
    throughout the pebble game.

    Also, the invariant:
    \begin{enumerate}
        \item {\bf (I4)}
        $\grsp_R V'+\out_R V'+\peb V'\le \grsp_R V'+\out_H V'+\peb V'$
    \end{enumerate}
    holds for all $V'\subset V$.
\end{lemma}
\begin{proof}
The proof of {\bf (I1)}, {\bf (I2)} and {\bf (I3)} are similar to the
proof of \reflem{pebble-game-invariants}.

For {\bf (I4)}, we just need to observe that since $E_H(V')\subset E_R(V')$,
the out degree in $H$ it at least the out-degree in $R$.
\end{proof}

From \reflem{represented-pebble-game-invariants} we see that $R$ must be
sparse, and by construction $R$ has dimension at least $(\ell+1)/k$.
Since $R$ is a pebble game graph, we see that $G$ is critical if and
only if $G=R$ for every represented pebble game construction.

\begin{restate}{representation}{{\bf Critical Representations}}
     $G$ is a critical sparse hypergraph of
    dimension $s$  if and only if the representation found by the pebble
    game construction coincides with $G$. This implies that $G$ is
    $s$-uniform and $\ell=sk-1$.
\end{restate}
\begin{proof}
The theorem follows from the fact that we can always move pebbles between
the ends of an independent set of vertices unless there are exactly $sk$ pebbles on it already, which is exactly the acceptance condition for the $(k,sk-1)$-pebble game.
\end{proof}

The observation that $E_H(V')\subset E_R(V')$ also proves that
any component in $H$ induces a block in $R$.
It is instructive to note that blocks in $R$ do {\it not} necessarily
correspond to blocks in $H$.

\section{Conclusions and Open Questions}
We have generalized most of the known results on sparse graphs to the
domain of hypergraphs.  In particular, we have provided graph theoretic, algorithmic
and matroid characterizations of the entire family of sparse hypergraphs
for $0\le \ell<ks$.

We also provide an initial result on the meaning of dimension in sparse hypergraphs;
in particular the representation theorem shows that the sparse hypergraphs for
$l\ge 2k$ are somehow intrinsically not $2$-dimensional.

The results in this paper suggest a number of open questions, which we
consider below.

\paragraph{Algorithms.} The running time and space complexity of the pebble
game with components is the natural generalization of the $O(n^2)$
achieved by Lee and Streinu in  \cite{LeSt05}.  Improving our
$\Omega(n^{s^*})$ running time to $O(m+n^2)$ may be possible with a
better data structure.

For the case where $d=2$, the pebble games of Lee and Streinu are not the
best known algorithms for the maps-and-trees range of parameters.  We do
not know if the algorithms of  \cite{GaWe88} and   \cite{gabow1995}
generalize easily to hypergraphs.

\paragraph{Graph theory.}

Proving a partial converse of  the lower-dimensional representation
theorem \refthm{representation}
is of particular interest to a
number of applications in rigidity theory.

%%%%%%
\bibliographystyle{abbrv}
\bibliography{pg_arboricity}

\appendix

\begin{center}
{\bf \Large Appendix}
\end{center}

\section{The matroid of sparse hypergraphs}
In this section we investigate matroidal properties of the sparse
graphs.  The main result of this section is due to  White and
Whiteley \cite{whiteley:matroids} where it is proven using the
circuit axioms. For completeness, we include another proof using the
basis axioms.

\begin{theorem}\labelthm{kl-matroid}
Let $\mathcal{B}$ be the collection of all tight graphs on $n$
vertices.  Then $\mathcal{B}$ is not empty when $k$, $\ell$, $n$ and
$d$ meet the conditions of \refthm{good-range} and $\mathcal{B}$ is
class of bases of a matroid $\mathcal{M}_{k,\ell}$ which has the
sparse graphs as its independent sets and the circuits as described
in \refsec{preliminaries} as its circuits.
\end{theorem}
\begin{proof}%[Proof of \refthm{kl-matroid}]
We verify that $\mathcal{B}$ obeys the basis axioms.  For
completeness, we state them here.

\begin{description}
\item[{\bf (B1)}] $\mathcal{B}\neq \emptyset$
\item[{\bf (B2)}] All bases are have the same cardinality.
\item[{\bf (B3)}] For distinct bases $B_1$ and $B_2$ there are
elements $e_1\in B_1-B_2$ and $e_2\in B_2-B_1$ such that
$B_1-e_1+e_2$ is a base.
\end{description}

{\bf (B1)} Follows from \refthm{good-range}.

{\bf (B2)} All tight graphs have exactly $kn-\ell$ edges.

{\bf (B3)}  Let $B_1$ and $B_2$ be distinct bases.  Then $B_1-B_2$
is not empty; let $e_2$ be an element of $B_1-B_2$ or dimension $s$.
Let $C$ be the subgraph induced by the vertex intersection of every
block in $B_{1}$ spanning $e_{2}$; $C$ is well-defined since $B_{1}$
is a block, and by \refthm{block-structure}, $C$ is a block.  (In particular,
$C$ is the inclusion-wise minimal block containing $e_2$.)
Moreover, $C-e_{2}$ is not empty; by hypothesis $C$ cannot be
$sk-\ell$ copies of $e_{2}$.

A graph that contains a subgraph that is not sparse called
{\bf dependent}.
Observe that any dependent subgraph in
$B_{1}+{e_{2}}$ must contain $C+e_{2}$.  By construction, no
subgraph of $C$ is tight, and thus $e_{2}$ is independent of any
subgraph of $B_{1}$ not containing $C$.

Let $e_{1}$ be an edge in $C-e_{2}$.  By the previous observation,
$C-e_{1}+e_{2}$, and thus $B_{1}-e_{1}+e_{2}$ is sparse.
\end{proof}
\end{document}